\documentclass[11 pt]{article}
\usepackage{sibarticle_SY}

\usepackage{ulem}

\usepackage{listings}
\usepackage{caption}
\usepackage{algorithm}
\usepackage{algpseudocode}
\usepackage{mathrsfs}
\usepackage{cancel, breqn, diagbox}
\newcommand{\commOut}[1]{}

\newcommand{\setR}{\mathbb{R}}

\newcommand{\expec}{\mathbb{E}}

\newcommand{\momentumParam}{\gamma}

\setuptodonotes{backgroundcolor=red!10, linecolor=red,
  bordercolor=red, size=\tiny, shadow, tickmarkheight=0.1cm}

\newcommand{\bm}[1]{\mbox{\boldmath{$#1$}}}

\newcommand{\CD}{\mathcal{D}}

\newcommand{\CL}{\mathcal{L}}
\newcommand{\CP}{\mathcal{P}}

\newcommand{\vb}{\bm{b}}
\newcommand{\vc}{\bm{c}}

\newcommand{\vg}{\bm{g}}
\newcommand{\vp}{\bm{p}}

\newcommand{\vs}{\bm{s}}

\newcommand{\vu}{\bm{u}}

\newcommand{\vx}{\bm{x}}

\newcommand{\zv}{\bm{0}}

\newcommand{\valpha}{\bm{\alpha}}
\newcommand{\vbeta}{\bm{\beta}}

\newcommand{\vlambda}{\bm{\lambda}}

\newcommand{\mA}{\bm{A}}
\newcommand{\mB}{\bm{B}}

\newcommand{\tr}{^{\intercal}}

\newcommand{\hsp}{\hspace{4mm}}

\newcommand{\hsppp}{\hspace{12mm}}

\newcommand{\maximize}{\mbox{maximize\hspace{4mm} }}
\newcommand{\minimize}{\mbox{minimize\hspace{4mm} }}
\newcommand{\subto}{\mbox{subject to\hspace{4mm}}}

\newcounter{commentcounter}
\long\def\symbolfootnote[#1]#2{\begingroup
  \def\thefootnote{\fnsymbol{footnote}}\footnote[#1]{#2}\endgroup}

\newcommand{\propnum}[2]{\vspace{3mm}
  \noindent {\bf Proposition #1}{\it #2} \vspace{3mm}}

\newcommand{\thmnum}[2]{\vspace{3mm}
  \noindent {\bf Theorem #1}{\it #2} \vspace{3mm}}
  \newcommand{\cornum}[2]{\vspace{3mm}
  \noindent {\bf Corollary #1}{\it #2} \vspace{3mm}}

\title{Differentially Private Linear Optimization for Multi-Party Resource Sharing} 
\ShortTitle{Differentially Private Resource Sharing}
\ShortAuthors{Karaca, Aydın, Yıldırım \& Birbil }

\NumberOfAuthors{4}
\FirstAuthor{Utku Karaca}
\FirstAuthorAddress{Erasmus University Rotterdam, 3000 DR, Rotterdam P.O. Box 1738, The Netherlands}
\FourthAuthor{Ş. İlker Birbil}
\FourthAuthorAddress{University of Amsterdam, 11018 TV, Amsterdam P.O. Box 15953, The Netherlands}

\SecondAuthor{Nurşen Aydın}
\SecondAuthorAddress{University of Warwick, Coventry CV4 7AL, United Kingdom}

\ThirdAuthor{Sinan Yıldırım}
\ThirdAuthorAddress{Sabancı University, Istanbul, Turkey}


\allowdisplaybreaks
\begin{document}

\maketitle
{
\singlespace
\begin{abstract}
This study examines a resource-sharing problem involving multiple parties that agree to use a set of capacities together. We start with modeling the whole problem as a mathematical program, where all parties are required to exchange information to obtain the optimal objective function value. This information bears private data from each party in terms of coefficients used in the mathematical program. Moreover, the parties also consider the individual optimal solutions as private. In this setting, the concern for the parties is the privacy of their data and their optimal allocations. We propose a two-step approach to meet the privacy requirements of the parties. In the first step, we obtain a reformulated model that is amenable to a decomposition scheme. Although this scheme eliminates almost all data exchanges, it does not provide a formal privacy guarantee. In the second step, we provide this guarantee with a locally differentially private algorithm, which does not need a trusted aggregator, at the expense of deviating slightly from the optimality. We provide bounds on this deviation and discuss the consequences of these theoretical results. We also propose a novel modification to increase the efficiency of the algorithm in terms of reducing the theoretical optimality gap. The study ends with a numerical experiment on a planning problem that demonstrates an application of the proposed approach. As we work with a general linear optimization model, our analysis and discussion can be used in different application areas including production planning, logistics, and revenue management. 
\end{abstract}
}

\section{Introduction} \label{sec: Introduction}
Efficient use of resources is one of the major objectives in every industry. Through collaboration, companies can coordinate their activities, increase their utilizations, and obtain significant savings. For instance, in freight logistics, companies share less-than-full vehicles to compound the problem of excess capacity \citep{speranza2018trends}. This collaboration helps them to reduce transportation costs and improve their operational efficiencies. Similarly, airlines provide joint services to increase their flight capacity utilization \citep{topaloglu2012duality}. In the automotive industry, Toyota and Fuji Heavy Industries shared an existing manufacturing site to increase the efficiency in manufacturing \citep{toyota06}. In supply chains, collaborative planning improves production and replenishment decisions \citep{poundarikapuram04}. There are also other examples in lodging, telecommunication, and maritime transportation where collaborations among companies occur by sharing resources, exchanging information, or providing joint services \citep{guo2018capacity, agarwal2010network, chun2016when}.

Although coordinating activities and sharing resources can provide significant benefits to companies, it also poses several challenges. These partnerships are mainly built around information exchange to coordinate the collective decision-making process. However, individual partners, though often working towards a common goal, can be competitors and may be unwilling or unable to fully disclose sensitive information about their operations \citep{hyndman2013aligning, Albrecht15}. Sensitive information may involve demand forecasts, selling prices, operational costs, and available capacities. In air-cargo transportation, for example, airlines and freight forwarders collaborate to sell the flight capacity. In this partnership, the freight forwarders keep their demand information, operating costs, and reservation prices private to protect their interests \citep{amaruchkul11}. While centralized decision-making is desirable to achieve maximum gain from the collaboration, it may not be realistic for collaborative planning due to restrictions in information exchange. Therefore, decomposition and decentralization approaches have been studied to minimize information sharing among participants \citep{Albrecht15, Ding19}. In decentralized systems, individual parties release only the required information at each step to obtain the optimal allocations of the shared resources with an iterative solution approach. This information can also be exchanged via a central planner. Although information exchange is reduced with a decentralization approach, the parties still have to share information about their operations with others or the central planner. Specifically, in a network resource-sharing setting, this shared information can be related to the optimal capacity allocations or the bid (dual) prices, which may reveal private information, like individual profits and capacity utilizations \citep{poundarikapuram04, Albrecht15}. How to protect the privacy of such information becomes a critical issue in practice. While numerous existing privacy protection methods depend on transformation and encryption-based approaches, two main issues arise with these methods: either they fail to strictly ensure privacy or they become computationally infeasible when applied to even small-scale problems \citep{birbil2020data, Chen22}. The data privacy concern in collaborations raises an important question: How can one mathematically guarantee data privacy while conducting optimization in a multi-party resource-sharing setting?

\subsection{Contributions} \label{sec: Contributions}
In this work, we address the aforementioned issues and the question above for a multi-party resource-sharing model in a general linear optimization framework. We discuss next the main methodological and literature contributions of the paper.

We consider a general linear optimization framework, where parties make use of the shared and private resources to manage their operations. Aside from the shared resources, the parties do not want to reveal directly their private information consisting of the coefficients in the objective function and the individual constraints. To this end, we first introduce a decomposition approach coupled with an iterative solution algorithm where each party is only required to share their allotments at each iteration. Then, we mask each allocation by adding a specifically adjusted noise to randomize the output. This randomization allows us to give a data privacy guarantee by using the mathematically rigorous definition of \textit{differential privacy} \citep{dwork2006calibrating}. Differential privacy is a widely studied concept in the computer science literature and has become an established privacy standard in various industries. The key idea is that even if an adversary knows the output of a differentially private algorithm, making accurate inferences about the data of individuals used in the algorithm remains challenging. As we have an iterative algorithm, this theory dictates using only a finite number of iterations to preserve privacy at a given level. Considering the random output and the limited number of iterations, we observe that our differentially private decomposition algorithm may end up with a sub-optimal solution. Therefore, we provide theoretical bounds on the difference between the approximate and the optimal objective function values. 

The private resource-sharing methods proposed in this paper satisfy a particular version of differential privacy, namely \textit{local differential privacy}. It is a stronger version of the standard differential privacy which is typically satisfied with a trusted aggregator. In contrast, locally differentially private algorithms, like the ones studied in this paper, operate without the need for a trusted aggregator. To the best of our knowledge, this is the first formal treatment of data privacy \textit{without any trusted party} in multi-party resource sharing via linear optimization. Using differential privacy without the assumption of a trusted aggregator results in larger noise in the shared quantities. To circumvent that problem, we further propose a novel modification to increase the efficiency of the private decomposition algorithm by reducing the variance of the noise, and consequently, the optimality gap. 

We support our methods and theoretical results with a set of numerical experiments on a production planning example and present the trade-off between optimality and privacy guarantee. We specify our framework for linear programming as it is arguably the most frequently used optimization tool in practice \citep{ibaraki1988resource, topaloglu2012duality, Albrecht15}. We believe that our analyses could inspire further studies on privacy-aware operations management.

\subsection{Review of Related Literature} \label{sec: Review of Related Literature}
Collaborative relationships between organizations have received considerable attention in the literature, and it is very popular in many industries including airlines \citep{wright2010dynamic}, logistics, maritime shipping \citep{agarwal2010network, verdonck2013collaborative, Lai2019}, and retail \citep{guo2018capacity}. Most of these studies assume complete information exchange among partners, which may not be realistic in practice due to the regulations. Therefore, recent studies focus on the techniques to minimize information sharing in collaborative optimization problems. \cite{poundarikapuram04} propose a decentralized decision-making framework based on the L-shaped method for a collaborative planning problem in a supply chain. The centralized collaborative planning problem is separated into a master problem that includes the common variables for all parties and several sub-problems that include private local objectives and variables. The authors present an iterative procedure to address this problem, where the parties can solve their local problems privately and disclose limited information to solve the master problem at each iteration. \cite{topaloglu2012duality} focuses on capacity sharing in airline alliances and proposes a decomposition approach to find booking limits for each alliance partner as well as bid prices for shared flights. \cite{kovacs13} study collaboration in lot-sizing problems with two parties and investigate different solution approaches, including centralized and decentralized methods. To minimize information exchange, the lot-sizing problem is decomposed and solved by each party sequentially in the decentralization approach. \cite{Albrecht15} propose a scheme to coordinate the collaborative parties in a decentralized environment where parties' local private problems can be modeled as linear programs. In the proposed scheme, each party exchanges proposals by sharing information on the primal optimal solution of the linear programs. 

Decomposition approaches have also attracted great attention in the machine-learning community. Federated learning is initially proposed in \cite{shokri2015privacy, konevcny2016federated,mcmahan2017communication}, and it focuses on learning tasks in a distributed environment. Commonly, machine learning models require a central training phase where all the data needs to be shared. The use of decomposition approaches enables training locally so that there is no need for data exchange. Even if the data is not shared, adversarial attacks are still possible through the shared model parameters during the training process. 
There are many studies in the literature to provide privacy guarantees by employing various methods. \cite{geyer2017differentially} design a differentially private federated learning algorithm. \cite{truex2019hybrid, yang2019federated} study privacy of federated learning by employing both differential privacy and secure multi-party computation protocols. We refer the reader to \cite{li2020federated} and \cite{mothukuri2021survey} for extensive literature reviews on federated learning and its privacy, respectively.

One well-known approach to overcoming privacy breaches in data sharing through collaborations is using transformation-based methods. These methods allow the reconstruction of the given model with altered data while preserving the optimality. For instance, \cite{vaidya2009privacy} proposes a transformation method in linear programming models where one party owns the objective function, and the other owns the constraints. In follow-up work, \cite{bednarz2009hiccups} show that the proposed method by \cite{vaidya2009privacy} is valid under specific restrictions. Later, \cite{hong2014inference} have revised the initial approach of \cite{vaidya2009privacy} and improved the level of privacy. \cite{mangasarian2011privacy} and \cite{mangasarian2012privacy} present a transformation-based approach to preserve privacy in both vertically and horizontally partitioned linear models. In another study, \cite{dreier2011practical} show the computational inefficiency of the cryptographic methods in mathematical models and propose a transformation-based approach. Additionally, they provide a security analysis of their approach. All of these studies focus on hiding both input data and the primal decision variables in various settings. \cite{Hastings22} work on secure multi-party computation for financial networks and analyze the efficiency and convergence of their proposed algorithms. Recently, \cite{birbil2020data} study data-privacy in a collaborative network revenue management problem using a transformation-based approach. The proposed method keeps both primal and dual variables private to each party. The transformation-based methods have a flaw in that either they do not guarantee privacy, or they are computationally not feasible for solving even small-scale problems \citep{sweeney1997weaving, narayanan2006break, toft09}.
As a remedy, the notion of differential privacy, which was proposed in the computer science field \citep{Dwork_2006}, has been widely studied in the literature.

Differential privacy is a probabilistic definition of privacy that suggests using random perturbations to minimize the probability of revealing specific records from neighboring databases \citep{dwork2006calibrating}. Several methods with differential privacy guarantees have been adopted in industry  \citep{Apple17, Erlingsson14}. \cite{Chaudhuri11} study the differentially private empirical risk minimization in machine learning where the loss function and the regularizer are convex. Following the work of \cite{dwork2006calibrating}, they propose two algorithms satisfying differential privacy based on output and objective function perturbations, respectively. \cite{Raef14} consider special cases in convex empirical risk minimization with (non)-smooth regularizer. Extending the work of \cite{Chaudhuri11}, they propose a gradient perturbation approach to satisfy differential privacy. \cite{Wang18} also study the gradient perturbation approach for (non)-convex loss functions. \cite{zhang2017dynamic} present two differentially private algorithms using the alternating direction method of multipliers for the empirical risk minimization problem. Recently, differential privacy has gained significant attention in the community of operations research and management science. \cite{huang2015differentially} consider a cost minimization problem among parties with an additional privacy requirement for the cost function. They employ differential privacy and present a private distributed optimization scheme. \cite{hale2015differentially} study a multi-party optimization problem in a cloud structure, where the aim is to optimize global objective function with global inequality constraints. In the study, the parties send all the necessary data to the cloud, whose role is to keep each party's data private. \cite{han2016differentially} study a differentially private distributed constrained optimization problem in resource allocation where they assume a convex and Lipschitz objective function. They consider only shared resources in a specific structure, called the server-client network structure. Even though this study is conducted for a resource allocation setting, it differs from our study by the direct application of the stochastic subgradient method and its lack of rigorous bounds on approximate optimality. We also note that there is a growing body of literature on data privacy in service systems. \cite{Chen23}, \cite{Chen22} and \cite{Lei2023} study private sequential decision-making and dynamic pricing. They develop pricing algorithms with differential privacy guarantees for online and offline settings.

\cite{hsu2014privately} study differential privacy in linear programs. They define variations of neighboring datasets and analyze solution methods for each of the definitions. Later, \cite{hsu2016jointly} focus on privacy in distributed convex programming models where the solution is naturally partitioned among the parties defining the problem. For this kind of problem, joint differential privacy (a {\it{relaxed}} definition of differential privacy) has been considered. It allows the output given to all other parties excluding one party to be insensitive to that party's input data \citep{Kearns14}. Joint differential privacy is usually achieved with the presence of a trusted aggregator. \cite{hsu2016jointly} study privacy in resource allocation problems using joint differential privacy with a trusted aggregator who collects data from all parties. This notion is widely used in mechanism design and equilibrium computation in certain large games where the aggregator (or mediator) collects data from players and recommends actions to form an equilibrium \citep{Cummings15, Kearns2015, Kearns14, Rogers14}.  Different from these studies, we consider the stronger notion of local differential privacy. In particular, we do not assume the existence of a trusted aggregator. The consequence of using local differential privacy --without the assumption of a trusted aggregator-- is larger noise in the shared quantities. We tackle this issue with smart step-size selection techniques as well as an adaptive clipping method that aims to reduce the noise while still preserving the same level of differential privacy.

\section{Methodology} \label{sec: Methodology}
In this section, we introduce the mathematical model for multi-party resource-sharing and our initial attempt to obtain a \textit{data-hiding} solution approach with decomposition. Suppose we have $K > 1$ parties that agree to collectively use $m \geq 1$ \textit{shared resources} with capacities given in the vector $\vc \in [0, \infty)^{m}$. In addition, each party $k$ has its $m_{k}$ \textit{private constraints}; such as individual resources, demand restrictions, production requirements, and flow conservation. We let $\vb_k \in \setR^{m_k}$ be the right-hand-sides of the private constraints of party $k$. Next, each party $k$ has matrices $\mA_k\in \setR^{m\times n_k}$ and $\mB_k\in \setR^{m_k\times n_k}$ that include the technology coefficients related to the shared and private constraints, respectively. Here, $n_{k}$ represents the number of variables controlled by party $k$, \textit{e.g.}, the number of products. The objective of the model is to find the optimal allocations of the parties such that the total utility is maximized (without loss of generality, we give a maximization model, and all our discussion can be trivially extended to a minimization setting). For each $k$, we let $\vx_{k} \in \mathbb{R}^{n_{k}}$ be the vector of party $k$'s individual variables and $\vu_{k} \in \mathbb{R}^{n_{k}}$ be its utility vector. The canonical multi-party resource-sharing model then becomes a linear program:
\begin{subequations}
\begin{eqnarray}
Z_P ~  = & \maximize  & \sum_{k = 1}^{K} \vu_k\tr \vx_k,     \label{eqn:org_model_obj} \\
	& \subto  & \sum_{k = 1}^{K} \mA_k\vx_k \leq \vc,        \label{eqn:org_model_c} \\
	 & & \mB_k\vx_k \leq \vb_k,  \quad k  = 1, \ldots, K. \label{eqn:org_model_nonneg}
\end{eqnarray}
\end{subequations}
The constraints in \eqref{eqn:org_model_c} guarantee that the total consumptions of the parties do not exceed the shared capacities whereas those in \eqref{eqn:org_model_nonneg} express the private constraints of each party $k$. We note that the linear model in \eqref{eqn:org_model_obj}-\eqref{eqn:org_model_nonneg} is already quite general and used in many real-world applications from production planning to portfolio selection, from logistics to network revenue management \citep{ibaraki1988resource, topaloglu2012duality, birbil2014network, Albrecht15}. 

In its current form, the model \eqref{eqn:org_model_obj}-\eqref{eqn:org_model_nonneg} can be solved when all input is available. Before raising the privacy concerns, let us first define formally what constitutes the (private) dataset for each party.
\begin{definition}[Dataset] \label{defn:dataset}
The dataset of party $k$ is defined as $\mathcal{D}_k = \{ \mA_{k}, \mB_{k}, \vb_{k}, \vu_{k} \}$.
\end{definition}
Although the parties agree to solve the resource-sharing problem collectively, they are unwilling to share their datasets with their fellow parties. This is because a dataset consists of sensitive information for its corresponding party. For instance, in a network revenue management problem, a dataset of a party may consist of private information about its products, capacities, and revenues \citep{birbil2020data}. Given this privacy concern, the main question becomes: How can a party join the others to solve the resource-sharing problem collectively, but keep its dataset private? To address this question, we first propose a decomposition approach, where each party ends up solving a sequence of sub-problems without sharing any data with the other parties.

\subsection{Data-hiding with Decomposition}  \label{sec: Data-hiding with Decomposition}
Our decomposition approach stems from a natural observation: If the shared capacities are divided among the parties \textit{a priori}, then there is no need to solve the overall problem as each party can independently solve its problem with the allocated capacity. However, ensuring a feasible allocation of shared resources is not always possible unless we have sensitive information from all parties. In addition, a priori partitioning of the shared capacities can easily lead to a sub-optimal solution without maximizing the total utility; see \textit{e.g.}, the discussion given by \cite{birbil2014network} on network capacity fragmentation. Therefore, we need to devise a scheme that divides the shared capacity optimally among the parties without revealing their private data.

Define the new decision vector $\vs_k \in \setR^{m}$ to be the allocation to party $k$ from the shared resources, and let $\vc = \sum_{k =1}^{K} \vs_k$. We can then write problem \eqref{eqn:org_model_obj}-\eqref{eqn:org_model_nonneg} equivalently as
\begin{subequations}
\begin{eqnarray}
	Z_\CP  ~= & \maximize & \sum_{k = 1}^{K} \vu_k\tr \vx_k, \label{eqn:sec_model_obj}     \\
	& \subto    & \mA_k\vx_k \leq \vs_k,   \quad k  = 1, \ldots, K,  \label{eqn:sec_model_Ak} \\
	&           &\mB_k\vx_k \leq \vb_k, \quad k  = 1, \ldots, K, \label{eqn:sec_model_nonnegprev} \\
	& 	      & \sum_{k = 1}^{K} \vs_k = \vc,  \label{eqn:sec_model_sk_c} \\
	&           & \vs_k \geq \zv, \quad  k  = 1, \ldots, K.  \label{eqn:sec_model_nonneg}
\end{eqnarray}
\end{subequations}
This model is almost separable except for the set of constraints \eqref{eqn:sec_model_sk_c}. We can relax this constraint by introducing a Lagrange multiplier vector $\vlambda\in\setR^{m}$. Then, the objective function of the relaxed problem becomes
\[
\CL(\vx, \vs, \vlambda) := \vc\tr \vlambda + \sum_{k = 1}^{K} \left( \vu_k\tr \vx_k - \vs_k\tr\vlambda \right),
\]
where $\vx := (\vx_1,\dots,\vx_K)$ and $\vs := (\vs_1,\dots,\vs_K)$. If for each party $k$ we further define the \textit{sub-problem}
\begin{equation} 
\begin{aligned}
	g(\vlambda; \CD_k) ~ := ~ & \maximize & \vu_k\tr \vx_k - \vs_k\tr\vlambda, \\
	& \subto    &\mA_k\vx_k \leq \vs_k,\\
	&           & \mB_k\vx_k \leq \vb_k,  \\
	&           &\vs_k \geq \zv,
\end{aligned} \label{eqn:subproblem}
\end{equation}
and assume that its dual has a feasible solution, then we obtain the Lagrangian dual problem
\begin{equation}
Z_\CL := \min\limits_{\pmb{\lambda}} \max_{\pmb{x}, \pmb{s}}
\CL(\vx,\vs,\vlambda) = \min\limits_{\pmb{\lambda}} \left\{\vc\tr\pmb{\lambda} + \sum_{k =1 }^{K} g(\vlambda; \CD_k)\right\}.
\label{eq:dual}  
\end{equation}
Since we are dealing with linear programs, the strong duality trivially holds, and hence $Z_\CP = Z_\CL$. This implies that we can solve the dual problem \eqref{eq:dual} instead of the primal problem \eqref{eqn:org_model_obj}-\eqref{eqn:org_model_nonneg}. We denote the optimal solution to \eqref{eq:dual} by $\vx^*$, $\vs^*$ and $\vlambda^*$.

It is important to observe that the dual problem becomes separable over the parties when the Lagrange multiplier vector $\vlambda$ is fixed. Moreover, the well-known subgradient algorithm for solving the dual problem is based on updating the Lagrange multiplier vector iteratively \citep{bertsekas2015convex}, where the update at iteration $t$ takes the form
\begin{equation}
\label{eq:subup}
\vlambda^{(t+1)} = \vlambda^{(t)} - \nu^{(t)} \left(\vc - \sum_{k =1}^{K} \vs_k^{(t)}\right).
\end{equation}
 In \eqref{eq:subup}, $\nu^{(t)}$ is the step-size at iteration $t$. Assuming that norms of subgradients and the distance $\|\vlambda^{(0)} - \vlambda^*\|$ are bounded, the subgradient algorithm is guaranteed to decrease the sub-optimality bound at the rate of $t^{-0.5}$ for any choice of $\nu^{(t)}$ as suggested in \citet{boyd2003subgradient}.

\begin{figure}
	\centering
	\includegraphics[width=0.86\textwidth]{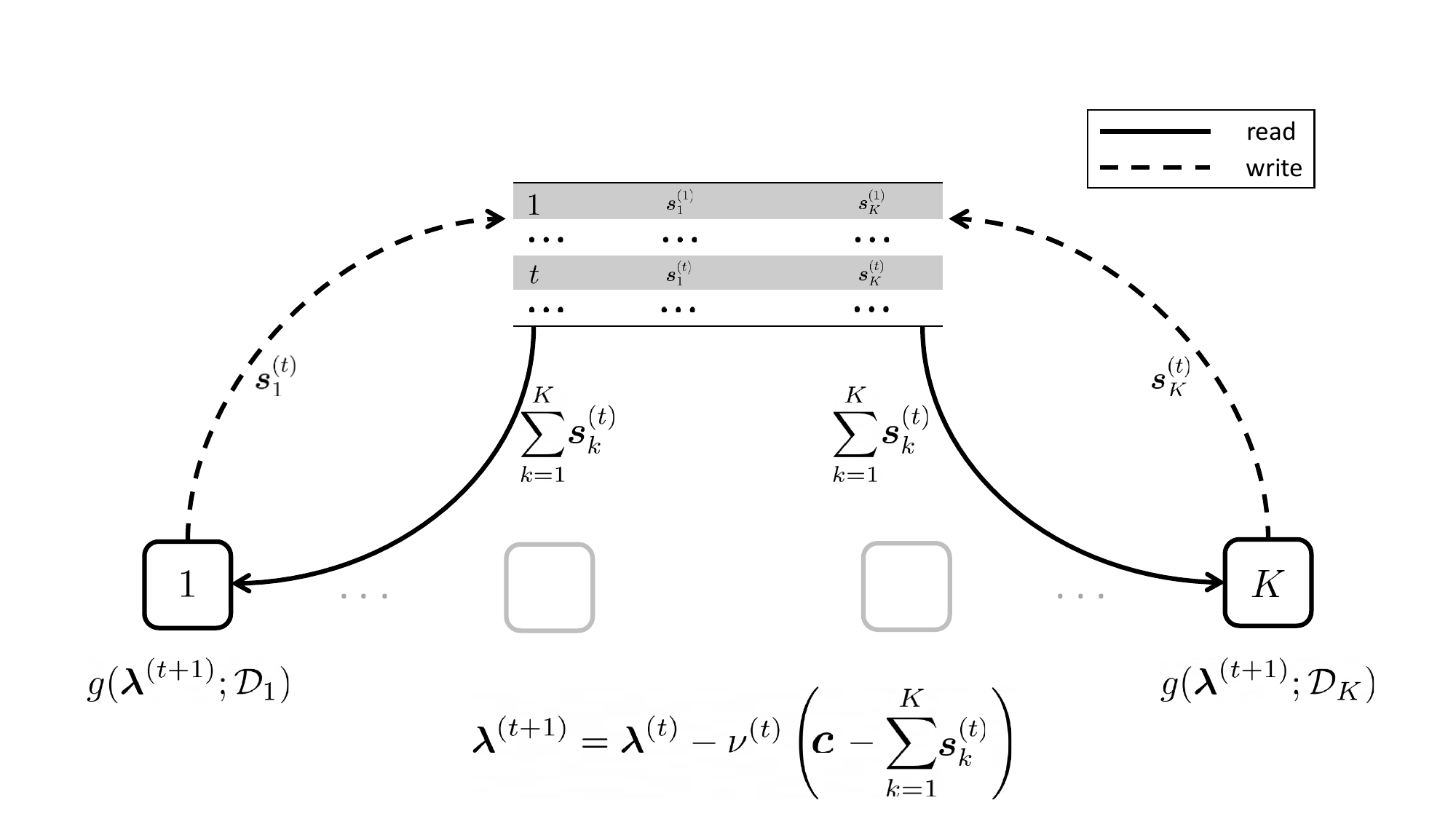}
	\captionof{figure}{Illustration of the dual decomposition approach to obtain the proposed data-hiding algorithm for multi-party resource sharing.}
	\label{fig:generalScheme}
\end{figure}

The outline of the dual decomposition approach is illustrated for one iteration in Figure \ref{fig:generalScheme}. The primal solutions $\vs_k^{(t)}$ are obtained by each party $k$ solving its sub-problem $g(\vlambda^{(t)}; \CD_k)$ on its own. Hence, there is no need for party $k$ to share its sensitive dataset $\CD_k$, but only the intermediate allocations $\vs_k^{(t)}$, with the other parties. These intermediate allocations are gathered to obtain $\vlambda^{(t+1)}$ as in \eqref{eq:subup}. The resulting update at iteration $t$ is given by
\begin{equation} \label{eq:subup0}
\begin{aligned}
(\vx^{(t)}_k, \vs^{(t)}_k) & = \arg \max g(\vlambda^{(t)}, \CD_k),   \hsp k = 1,\ldots,K, \\
\vlambda^{(t+1)} & = \vlambda^{(t)} - \nu^{(t)} \left(\vc - \sum_{k  = 1}^{K} \vs_k^{(t)}\right).
\end{aligned}
\end{equation}
The overall algorithm may start with the dual feasible solution $\vlambda = \zv$ (for ease of reference, we give the dual of \eqref{eqn:org_model_obj}-\eqref{eqn:org_model_nonneg} in Appendix \ref{sec:appDual}). 

In this setting, we assume that each party shares the obtained intermediate allocations with other parties truthfully. This is a well-known assumption when parties want to collaborate for their mutual benefit, and hence, act honestly to obtain the correct results \citep{Atallah03, hyndman2013aligning}. Furthermore, in many cases, the collaboration between firms lasts for more than one occasion. Therefore, a firm must consider the long-term prospects of future collaborations that can be jeopardized by opportunistic behavior. 

\paragraph{Updating with a momentum term:} 
One can consider adding a momentum term in the update step in \eqref{eq:subup} to improve the computational performance of the dual decomposition approach.
Accordingly, the new update step becomes
\begin{equation} \label{eqn:dp_momentum}
\vlambda^{(t+1)} = \vlambda^{(t)} - \nu^{(t)}\left(\vc - \sum_{k = 1}^{K} \vs_k^{(t)}\right) + \momentumParam \left(\vlambda^{(t)} - \vlambda^{(t-1)}\right).
\end{equation}
Here, $\momentumParam > 0$ is a momentum parameter. 


\subsection{Differentially Private Resource Allocation} \label{sec: Differentially Private Resource Allocation}

Although the data-hiding approach in Section \ref{sec: Data-hiding with Decomposition} prevents sharing sensitive data directly, it does not give a formal guarantee of data privacy. The deterministic and iterative nature of the approach could leak information regarding the datasets $\CD_{1}, \ldots, \CD_{K}$. To this end, one can consider applying transformation, encryption, or anonymization techniques to ensure data privacy. However, these methods either do not strictly guarantee privacy, or they are computationally not feasible for solving even small-scale problems \citep{sweeney1997weaving, narayanan2006break, toft09}. 

For a formal data privacy guarantee, we consider \textit{differential privacy} \citep{Dwork_2006}. In this section, we present a modified version of the data-hiding approach in Section \ref{sec: Data-hiding with Decomposition} that provides differential privacy. We adopt a particular version of differential privacy for our approach, namely \textit{local differential privacy}. We start by providing essential background material on local differential privacy and subsequently discuss its application to the multi-party resource allocation problem.

\subsubsection*{Background on Local Differential Privacy}

Local differential privacy (LDP) is a property that concerns randomized algorithms. A randomized algorithm takes an input dataset and returns a random output, where the randomness is intrinsic to the algorithm's inner mechanisms. Let $\mathscr{D}$ be the set of all possible datasets for a single party. A randomized algorithm that operates on datasets in $\mathscr{D}$ can be defined as a pair $\mathcal{A} = (\mu, A)$, where $\mu$ is a probability distribution on $\Omega$ representing the source of randomness in the algorithm, and the function $A: \Omega \times \mathscr{D} \mapsto \mathcal{O}$ returns the random output $A(\omega, \mathcal{D}) \in \mathcal{O}$ upon taking the random number $w \sim \mu$ and the dataset $\mathcal{D} \in \mathscr{D}$ as its inputs. Now, we give the formal definition of local differential privacy.

\begin{definition}[Local Differential Privacy] \label{defn:ldp}
Let $\epsilon, \delta > 0$. A randomized algorithm $\mathcal{A} = (\mu, A)$ is $(\epsilon, \delta)$-locally differential private if for all pairs of $\mathcal{D}, \mathcal{D}'$ and all $S \subseteq \mathcal{O}$, it satisfies
\[
\mathbb{P}[A(\omega, \mathcal{D}) \in S] \leq e^{\epsilon} \mathbb{P}[A(\omega, \mathcal{D}') \in S] + \delta,
\]
where the randomness is due to $\omega \sim \mu$.
\end{definition}
This definition implies an $(\epsilon, \delta)$-LDP algorithm outputs \textit{similarly} when the input dataset changes. Smaller $\epsilon$ and $\delta$ imply more privacy. 

Local differential privacy is a stronger version of the standard differential privacy, also referred to as centralized differential privacy. The difference between those two definitions is significant. Differential privacy concerns algorithms that operate on aggregate datasets from multiple users, such as $(\mathcal{D}_{1}, \ldots, \mathcal{D}_{k}, \ldots, \mathcal{D}_{K})$, and it ensures that the output does not change much when one user's data in the aggregate dataset changes, \textit{e.g.}, leading to a neighbor data collection $(\mathcal{D}_{1}, \ldots, \mathcal{D}_{k}', \ldots, \mathcal{D}_{K})$. That is why differential privacy, and some of its weaker variants such as joint differential privacy, typically require a trusted aggregator. Local differential privacy, on the other hand, concerns algorithms that operate on a single party's data. In loose terms, each party protects its privacy by sharing a noisy version of a computation performed solely on its own data. In this way, locally differentially private algorithms do not require a trusted aggregator. Apart from this difference, the same mathematical tools can be used for both types of differential privacy. In particular, properties such as composition and post-processing apply to either type. Likewise, noise-adding mechanisms, like the Gaussian mechanism, can be used to provide both types of privacy. Finally, relations between other forms of privacy, such as Theorem \ref{thm:zCDPtoLDP}, apply similarly to both types of differential privacy.

In the following, we provide the tools of (local) differential privacy that we use in our work. In particular, we define zero-concentrated differential privacy (zCDP) \citep{bun2016concentrated} and provide some results on its relation to LDP. 
We consider zCDP because it is well-suited to the Gaussian mechanism. With a slight abuse of terminology, we will adapt the following definitions and results for LDP. We will consider randomized algorithms that operate on a single user's dataset in $\mathscr{D}$. 
\begin{definition}[zCDP \citep{bun2016concentrated}] \label{defn:zCDP} 
$\mathcal{A}$ is $\rho$-zCDP if for all $\mathcal{D}, \mathcal{D}' \in \mathscr{D}$ and all $\alpha\in(1,\infty)$, we have
\[
D_{\alpha}(A(\omega, \mathcal{D})|| A(\omega, \mathcal{D}')) \leq \rho\alpha,
\]
where $\omega \sim \mu$ and $D_{\alpha}(A(\omega, \mathcal{D})|| A(\omega, \mathcal{D}'))$ is the $\alpha$-R\'{e}nyi divergence between the distributions of $A(\omega, \mathcal{D})$ and $A(\omega, \mathcal{D}')$. 
\end{definition}

The next result gives the relationship between (local) differential privacy and zero concentrated (local) differential privacy that is due to \cite{bun2016concentrated}.

\begin{theorem}[LDP and zCDP {\citep[Proposition 3]{bun2016concentrated}}]
\label{thm:zCDPtoLDP}
 If $\mathcal{A}$ is $\rho$-zCDP, it is $(\epsilon, \delta)$-LDP for all $\delta > 0$ and $\epsilon =\rho + \sqrt{4\rho\log(1/\delta)}$.
\end{theorem}

Most (locally) differentially private algorithms achieve randomization by adding noise to their output. We consider the Gaussian mechanism in this paper, where the output of a function of the sensitive data is distorted with a noise drawn from the normal distribution. The variance of the normal distribution should be adjusted according to the $L_{2}$ sensitivity of the function \citep{dwork2006calibrating}.
\begin{definition}[Sensitivity]
The $L_{2}$-sensitivity of a function $f: \mathscr{D} \mapsto \mathbb{R}^{d}$, $d \geq 1$, is given by
\[
\Delta_{f} = \displaystyle \sup_{ \mathcal{D}, \mathcal{D}' \in \mathscr{D}} \| f(\mathcal{D})-f(\mathcal{D}') \|_{2}.
\]
\end{definition}
\begin{theorem}[Gaussian Mechanism {\citep[Proposition 6]{bun2016concentrated}}] \label{thm:gausszCDP}
Suppose that $A(\omega, \mathcal{D}) = f(\mathcal{D}) + \omega$ and $\mu = \mathcal{N}(0, \sigma^{2} I)$. Then $\mathcal{A}$ is $\Delta_{f}^2/2\sigma^2$-zCDP.
\end{theorem}
When a differentially private algorithm outputs only a single value, verifying its privacy guarantee is straightforward. However, in practical scenarios, differentially private algorithms often entail multiple analyses, and it becomes necessary for the entire output to maintain differential privacy. In computer science literature, the process of combining multiple outputs from differentially private algorithms is termed composition. The following composition theorem quantifies the total privacy loss when multiple calculations of the sensitive data are released.
\begin{theorem}[Composition {\citep[Lemma 7]{bun2016concentrated}}] \label{thm:zCDPComp}
Let $\mathcal{A}_{i} = (\mu_{i}, A_{i})$, with $A_{i}: \Omega \times \mathscr{D} \mapsto \mathcal{O}_{i}$, is $\rho_i$-zCDP for $i = 1, \ldots, m$. Then, the composition $\mathcal{A} :=(\mu_{1} \otimes \ldots \otimes \mu_{m}, A: \mathscr{D} \mapsto \mathcal{O}_{1} \otimes \ldots \otimes \mathcal{O}_{m})$ where 
\[
A((\omega_{1}, \ldots, \omega_{m}), \mathcal{D}) = (A_{1}(\omega_{1}, \mathcal{D}), \ldots, A_{m}(\omega_{m}, \mathcal{D})), \quad \omega_{i} \sim \mu_{i}, \quad i = 1, \ldots, m
\]
is $\sum_{i = 1}^{m} \rho_i$\,-\,zCDP. Moreover, the same privacy guarantee holds even if $\mathcal{A}_{i}$ depends on the outputs of its predecessors, $\mathcal{A}_{1}, \ldots, \mathcal{A}_{i-1}$ for all $i = 1, \ldots, m$.
\end{theorem}

\subsubsection*{LDP Resource Allocation: Algorithmic Details}
Now we can present our LDP resource allocation algorithm. Consider the update in \eqref{eq:subup0}. As the parties collaborate, they share the data-sensitive outputs $\vs_{k}^{(t)}$ with each other. Hence, the resource allocation algorithm can be made LDP if each party $k$ uses the Gaussian mechanism to share their $\vs_k^{(t)}$'s with privacy-preserving noise. Specifically, the update in \eqref{eq:subup0} is modified as
\begin{subequations}
\begin{align} 
(\vx^{(t)}_k, \vs^{(t)}_k) & = \arg \max g(\vlambda^{(t)}, \CD_k),   \hsp k = 1,\ldots,K, \label{eq:subupdp1} \\
\tilde{\vs}_k^{(t)} &= \vs_k^{(t)} + \bm{\omega}_k^{(t)}, \quad k = 1, \ldots, K, \label{eq:subupdp}\\
\vlambda^{(t+1)} & = \vlambda^{(t)} - \nu^{(t)} \left(\vc - \sum_{k  = 1}^{K} \tilde{\vs}_k^{(t)}\right) + \momentumParam \left(\vlambda^{(t)} - \vlambda^{(t-1)}\right). \label{eq:subupdp2}
\end{align}
\end{subequations}

Algorithm \ref{alg:dpra} shows the steps of our differentially private algorithm. If we replace $\vs_k^{(t)}$ with $\tilde{\vs}_k^{(t)}$ in Figure \ref{fig:generalScheme}, then we obtain one iteration of the main loop in Algorithm \ref{alg:dpra}.

\begin{algorithm}
\caption{Locally Differentially Private Resource Allocation Algorithm}\label{alg:dpra}
\begin{algorithmic}
\Require $\vlambda^{(0)} = \zv$; number of iterations $T$; privacy parameters $(\epsilon, \delta)$; \\
\hsppp~~ datasets $\mathcal{D}_{k} = \{ \mA_{k}, \mB_{k}, \vb_{k}, \vu_{k} \}$, $k = 1, \ldots, K$; step sizes $\nu^{(t)}$, $t = 1, \ldots, T$
\Ensure $(\vx^{(T)}_k, \vs^{(T)}_k), \, k=1,\ldots,K $
\For{$t = 0,1,\dots,T$}
\State Each party $k$ solves \eqref{eqn:subproblem} and obtains its allocation, $\vs_k^{(t)}$ as in \eqref{eq:subupdp1}
\State Each party $k$ publishes noisy allocation, $\tilde{\vs}_k^{(t)}$ as in \eqref{eq:subupdp}
\State Lagrange multiplier, $\vlambda^{(t+1)}$ is calculated by all parties with \eqref{eq:subupdp2}
\EndFor
\end{algorithmic}
\end{algorithm}

To ensure that this algorithm is $(\epsilon, \delta)$-LDP, the noise distribution in \eqref{eq:subupdp} must be adjusted to the parameters $\epsilon, \delta$, $m$, $T$, and the sensitivity of the evaluations. We assume that each party $k$ has an additional set of constraints $\vs_k \leq \vc$ added to  its set of private constraints $\{ \vx_{k} \in \mathbb{R}^{n_{k}}: \mB_k\vx_k \leq \vb_k \}$. This is a reasonable assumption as the shared capacities is known by all the involved parties. We can therefore ensure that
\begin{equation} \label{eq: bound on s}
0 \leq s_{k, j} \leq c_{j}, \quad j = 1, \ldots, m, ~ k = 1, \ldots, K.
\end{equation}
Consequently, the $L_{2}$-sensitivity of the mechanism in \eqref{eq:subupdp1} that computes $s_{k, j}$ is $c_{j}$. Furthermore, if the updates in \eqref{eq:subupdp1}-\eqref{eq:subupdp2} are run for $T$ iterations, we end up with a composition of $T m$ Gaussian mechanisms. The $L_{2}$-sensitivity for each component, together with a total number of iterations $T$, determine the variance of $\bm{\omega}_{k}^{(t)}$ to provide  $(\epsilon, \delta)$-LDP. The following proposition specifies that variance.

\begin{proposition}
\label{prop:diffPriv}
The multi-party resource sharing algorithm using updates \eqref{eq:subupdp} for $T$ iterations provides $(\epsilon, \delta)$-LDP if the random vectors are drawn as 
\begin{equation} \label{eq: capacity-specific DP variances-1}
\bm{\omega}_{k, j}^{(t)} \sim \mathcal{N}\left(0, \frac{T m  c_{i}^{2}}{2 \rho(\epsilon, \delta)} \right), \quad k = 1, \ldots, K; ~ j = 1, \ldots, m; ~t=1, \dots, T,
\end{equation}
where $\rho(\epsilon, \delta) = (\sqrt{\log(1/\delta) + \epsilon} - \sqrt{\log(1/\delta)})^{2}$.
\end{proposition}
The proof of Proposition \ref{prop:diffPriv}, given in Appendix \ref{sec:appProof}, exploits the zCDP of the Gaussian mechanism, the composition property of zCDP, and finally a conversion from zCDP to standard $(\epsilon, \delta)$-LDP.

The differential privacy noise variance given in Proposition \ref{prop:diffPriv} is calculated according to the component-wise sensitivities. It is also possible to use a common noise variance \textit{for all} $s_{k, j}$ values by considering the $L_{2}$-sensitivity of the whole vector $\vs_{k}$. Consider now the function  $F_{k, \lambda}: \mathscr{D} \mapsto [0, \infty)^{m}$ resulting in
\begin{equation} \label{eq: overall function}
F_{k, \lambda^{(t)}}(\mathcal{D}_{k}) = (s_{k, 1}^{(t)}, \ldots, s_{k, m}^{(t)}).
\end{equation}
It can be shown that the $L_{2}$-sensitivity of $F_{k, \lambda}$ is $\| \vc \|_{2}$ for all $k = 1,\ldots, K$ and $\lambda$. Hence, noise variance in \eqref{eq: capacity-specific DP variances-1} can be replaced by $T  \|\bm{c}\|^{2} /( 2\rho(\epsilon, \delta))$. The change here is the replacement of $m c_{j}^{2}$ by $\| \vc \|_{2}^{2}$. In fact, this is not the only other choice for noise variances. One can optimize this choice, \textit{e.g.}, by minimizing the sum of the variances along the components under the constraint of providing a certain amount of privacy per iteration. The optimization can be done in the same spirit as in \cite{kuru2020differentially}, where a similar noise-distribution problem is addressed. We do not pursue this point further and continue with the choice specified in Proposition \ref{prop:diffPriv}.

\subsection{Error bounds for LDP Resource Allocation Algorithm}
The privacy guarantee of the modified algorithm comes at the expense of stopping earlier possibly without reaching the optimal solution. Moreover, the stochasticity due to using noisy subgradients also affects the convergence behavior of the algorithm. We give an upper bound implying the suboptimality of the differentially private algorithm in the following theorem. The theorem requires a typical assumption used in the analysis of subgradient methods \citep{boyd2003subgradient}.
\begin{assumption} \label{asmp: bound on lambda}
$\|\vlambda^{(0)} - \vlambda^*\|_{2} \leq M$ for some $M > 0$. 
\end{assumption}

We define the distance to the optimal objective function value at iteration $t$ as
\[
G^{(t)} := \CL(\vx^{(t)}, \vs^{(t)}, \vlambda^{(t)}) - \CL(\vx^*, \vs^*, \vlambda^*), \quad t \geq 0.
\]
The following theorem establishes an upper bound for the expectation of $G^{(t)}$ without momentum updates, \textit{i.e.}, $\gamma = 0$.

\begin{theorem} \label{thm:convergence}
Suppose Assumption \ref{asmp: bound on lambda} holds and Algorithm \ref{alg:dpra} is run for $T$ iterations with a fixed step-size $\nu^{(t)} = \nu$ and with noise variances given in \eqref{eq: capacity-specific DP variances-1} to provide $(\epsilon, \delta)$-LDP. Then, there exists a $\nu$, for which the minimum expected distance to the optimal objective function value is bounded as
\[
\min_{t=0, 1, \ldots, T}~ \expec[G^{(t)}] \leq \frac{M B}{\sqrt{T}},
\]
where the constant $B$ is defined as
\begin{equation} \label{eq: B bound}
B = \|\vc\| \left(\frac{T m K}{2\rho(\epsilon, \delta)}+(K-1)^2 \right)^{1/2}.
\end{equation}

\end{theorem}
The proof of this result is given in Appendix \ref{sec:appProof}. Unlike in many applications using differential privacy, we observe that the bound increases with $K$.  This is because the variance of the overall noise with $K$ parties is proportional to $K$, while the sum of the non-noisy $\vs_{k}^{(t)}$ values target $\vc$. 
That is why summing does not necessarily improve, unlike the standard case, where the $K$ variables are independent with a distribution independent of $K$. This highlights the unique challenges and considerations in our scenario, where the standard assumptions of differential privacy may not directly apply due to specific constraints and characteristics of our problem setting.

Note that $\rho(\epsilon, \delta) \approx \epsilon^{2} / (4 \log(1/\delta))$. Hence, the bound in Theorem \ref{thm:convergence} could be reduced by increasing $\epsilon$ or $\delta$, making the algorithm less private. This is a typical trade-off for differentially private iterative algorithms, where either the privacy loss or the error increases with the number of computations using the sensitive data \citep{ji2014differential}.



When the momentum term is involved in the update \eqref{eq:subupdp2}, \textit{i.e.}, $\gamma > 0$, it is also possible to give a more specific bound on the suboptimality of Algorithm \ref{alg:dpra}. The proof is deferred to Appendix \ref{sec:appProof}.

\begin{theorem} \label{thm:convergence_tight}
Suppose Assumption \ref{asmp: bound on lambda} holds. If the multi-party resource-sharing algorithm with momentum updates is run with a fixed step-size for $T$ iterations with $(\epsilon, \delta)$-LDP, the minimum expected distance to the optimal objective function value is bounded as
\[
\min_{t=0, 1, \ldots, T}~ \expec [G^{(t)}] \leq \frac{M^2}{2 T \nu} + \frac{\nu}{2} B^{2}
\]
provided that the step-size satisfies $0 \leq \nu \leq \frac{\sqrt{ G^{(0)2} + (1 - \gamma) B^{2}T M^{2}} - G^{(0)}}{B^{2} T}$ with $B$ as in \eqref{eq: B bound}.
\end{theorem}

\paragraph{Remark.}\textit{With its current form, our differentially private algorithm uses the solution obtained at the last iteration, $T$. This is because parties only share $\vs^{(t)}_{k}$ at any iteration $t$, and there is no way to assess whether the best solution is achieved in that iteration without breaching privacy. However, the best solution could be tracked by extending the current sharing scheme such that each party shares their objective function values in each iteration, $g(\vlambda^{(t)};\mathcal{D}_k)$ in addition to $\vs_{k}^{(t)}$. This additionally shared data-sensitive information would also consume the privacy budget and have an impact on the achieved results due to potential increase in the noise variances.}

\subsection{Differential Privacy with Clipping and Adaptive Sensitivity} 
The bounds presented in Theorem \ref{thm:convergence} and Theorem \ref{thm:convergence_tight} become loose when $K$, the number of parties, increases. The main reason is as follows. Recall that for each $k=1, \ldots, K$ the sensitivity of $\bm{s}_{k}$ is $\bm{c}$. This results in the total variance of the noisy sum $\sum_{k = 1}^{K} \tilde{\vs}_{k}^{(t)}$ being proportional to $K$. To overcome that problem, we propose an alternative update method that reduces the sensitivity, thus reducing the amount of noise in $\tilde{\vs}_{k}$ values, without introducing additional privacy loss. Firstly, we determine a parameter $\alpha \geq 1$ and let $\bar{\vc} = \alpha \vc$. For iteration $t$, we introduce a capping variable, $\bar{\vs}_k^{(t)}$, with the initialization $\bar{\vs}_{k}^{(0)} = \bar{\vc}/K$. Instead of the update step presented in \eqref{eq:subupdp}, we have the following update steps:
\begin{subequations}
\begin{align}
& (\vx^{(t)}_k, \vs^{(t)}_k)  = \arg \max g(\vlambda^{(t)}, \CD_k),   \hsp k = 1,\ldots,K,  \label{eq:subupdpClip-1} \\
&\tilde{s}_{k, j}^{(t)}= \min\{ \bar{s}_{k, j}^{(t)}, s_{k, j}^{(t)} \} + w_{k, j}, \quad w_{k, j} \sim \mathcal{N} \left(0, \frac{T m \bar{s}_{k, j}^{(t)2}}{2  \rho(\epsilon, \delta)} \right),\ \ j =1,\ldots, m, \quad k = 1, \ldots, K, \label{eq:subupdpClip-2}\\
& \vlambda^{(t+1)}  = \vlambda^{(t)} - \nu^{(t)} \left(\vc - \sum_{k = 1}^{K} \tilde{\vs}_k^{(t)}\right)+ \momentumParam \left(\vlambda^{(t)} - \vlambda^{(t-1)}\right), \label{eq:subupdpClip-3} \\
& \bar{s}_{k,  j}^{(t+1)} = \bar{c}_{j} \frac{\max \{ \min\{c_{j},  \tilde{s}_{k, j}^{(t)} \}, \tau \}}{\sum_{k' = 1}^{K} \max \{ \min\{c_{j},  \tilde{s}_{k', j}^{(t)} \}, \tau \}}, \quad  j = 1, \ldots, m, \quad k = 1, \ldots, K, \label{eq:subupdpClip-4}
\end{align}
\end{subequations}
where $\tau$ is a small positive number (serving to avoid $0$).
In this form of updates, even though the value of $\vs_k$ can change from $\zv$ to $\vc$, we can make the sensitivity smaller. This is because the deterministic part in \eqref{eq:subupdpClip-2} bounded above by $\bar{s}^{(t)}_{k, j}$. Since $\sum_{k = 1}^{K}\bar{\vs}^{(t)}_{k, j} = \alpha c_{j}$, the amount of noise in $\sum_{k = 1}^{K} \tilde{s}_{k, j}^{(t)}$ is proportional to $\alpha c_{j}$. Since  $\alpha < K$, this amount of noise is lower than the $K \vc$ of the previous version of the algorithm. Additionally, in \eqref{eq:subupdpClip-4} we adjust $\bar{s}_{k, j}$ values in a predictive manner so that the capping value $\bar{s}_{k, j}^{(t)}$ for a party $j$ is made larger for the next iteration if its actual $s_{k, j}$ value, estimated by $\tilde{s}_{k, j}^{(t)}$ is larger. For completeness, we state that the modified algorithm that uses clipping to reduce the amount of privacy-preserving noise is still $(\epsilon, \delta)$-LDP; a proof is given in Appendix \ref{sec:appProof}.

\begin{corollary} \label{cor:difPrivClip}
The multi-party resource sharing algorithm using updates \eqref{eq:subupdpClip-1}-\eqref{eq:subupdpClip-4} for $T$ iterations is $(\epsilon, \delta)$-LDP.
\end{corollary}
As a result of the clipping and truncation operations in \eqref{eq:subupdpClip-1}-\eqref{eq:subupdpClip-4}, the approximate step is no longer an unbiased estimate of the subgradient. Therefore, our convergence and optimality-bound results above do not immediately apply. However, we have observed a favorable impact of smaller variance due to clipping and truncation in our numerical experiments.

\section{Computational Study} \label{sec:simStudy}

In this section, we present an application with the generic model \eqref{eqn:org_model_obj}-\eqref{eqn:org_model_nonneg}. We simulate a production planning problem with synthetic data to evaluate the efficiency of the proposed methods and discuss the effects of the privacy parameters on the results. We start by explaining our simulation setup in detail and then continue by presenting our numerical results. 

\subsection{Setup}
We construct an example with multiple parties and five shared capacities ($m = 5$). The components of the shared capacity vector $\vc\in\setR^{m}$ are randomly sampled from the interval $[10,20]$. Each party $k \in \{1, \ldots, K\}$ has $r_k$ private capacities, and $n_k$ products to produce. We also introduce demand constraints for each product, and hence, party $k$ has $m_k = n_k + r_k$ private constraints designated by the right-hand-side vector $\vb_k \in \setR^{m_k}$. The parameters $r_k$ and $n_k$ are sampled from the sets $\{5, 6, \dots, 10\}$ and $\{10, 11, \dots, 20\}$, respectively. The private capacity limits are randomly sampled from the interval $[0.0, 10.0]$. To guarantee to obtain a feasible problem with respect to the demand constraints, we first solve the problem to optimality without those constraints, and then, determine the intervals for randomly sampling the product demands. The components of the technology matrices $\mA_k \in \setR^{m \times n_k}$ and $\mB_k\in\setR^{m_k \times n_k}$ are obtained randomly from the intervals $[0.0, 5.0]$ and $[0.0, 1.0]$, respectively. The components of the utility vector $\vu_k \in \setR^{n_k}$ of party $k$ are sampled from the interval $[50.0, 150.0]$. 

We test the differentially private algorithm across various parameter settings, exploring both small and large-scale problems with the number of parties varying between $K \in {5,10}$. Aligned with industrial applications, we set the privacy budget at $\epsilon \in \{0.1, 0.5, 2.0, 4.0\}$ while keeping $\delta$ fixed at 0.001, \citep{Apple23, Erlingsson14}. The value of $\epsilon$, where $\epsilon > 0$, serves as a determinant of the extent of privacy leakage, with smaller values signifying stronger protection of participants' privacy. Our computational results are reported over 100 simulation runs for each experiment. In all our experiments for the differentially private setting, we use the constant step-size scheme. The choice of step-length plays a crucial role in influencing the convergence behavior of our differentially private algorithm due to using noisy subgradients. In our numerical experiments, we identify instances where the differentially private algorithms exhibit slower convergence or even divergence. To ensure a fair comparison among different privacy-preserving algorithms, we present results from the best-performing 90\% of 100 replications based on the optimality gap in all experiments. It is important to note that the primary objective of our numerical experiments is to assess the influence of privacy parameter $\epsilon$ on different problem settings, rather than to compete with state-of-the-art non-private solution algorithms.

\subsection{Results} \label{sec:simResults}
\begin{figure}[h]
	\centering
	\includegraphics[width=0.75\linewidth]{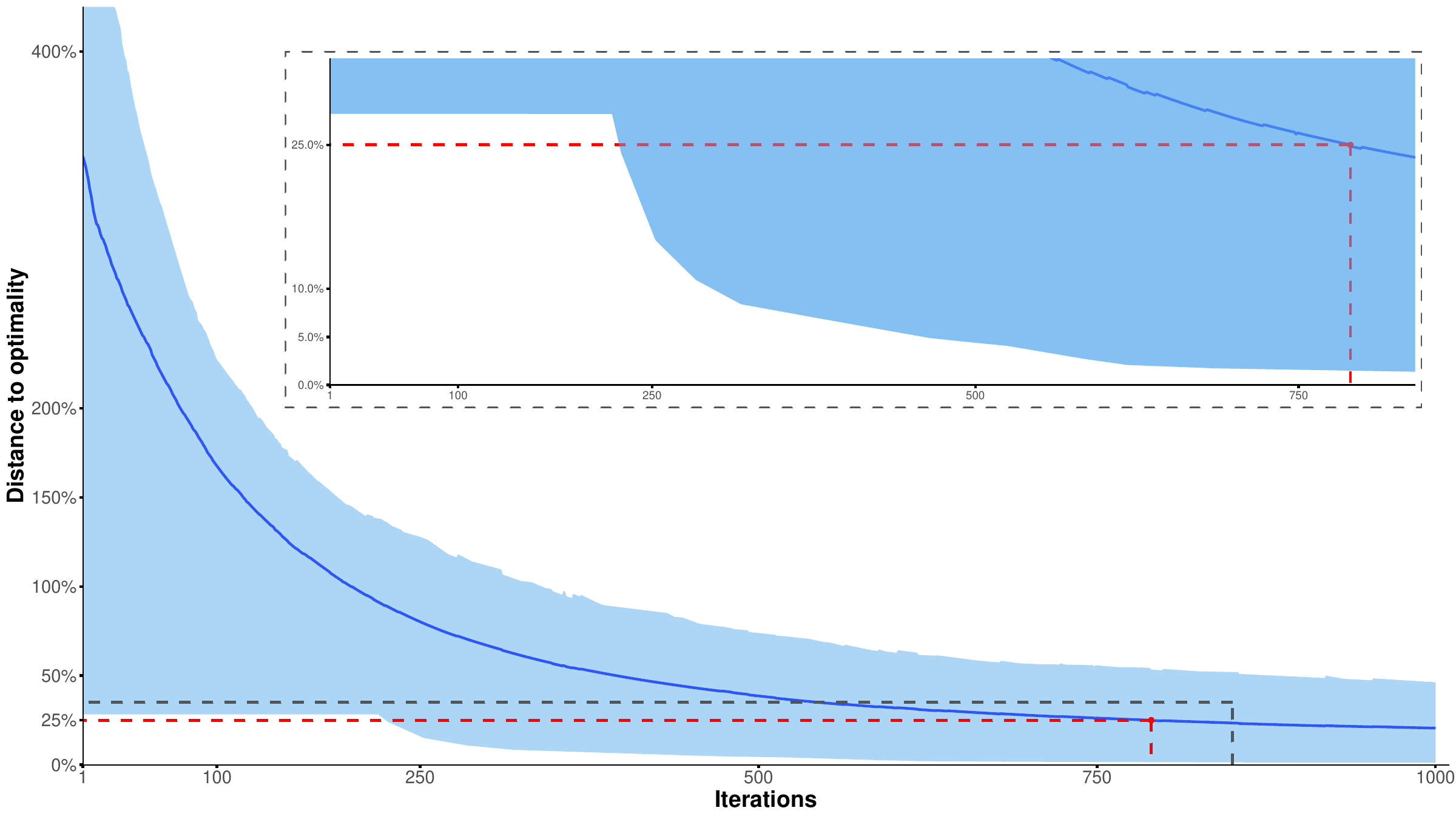}
	\caption{The percentage differences between the objective function values of the non-private model and the data-hiding model at different iterations. The subfigure shows an enlarged view of the rectangular region marked with the dashed black line on the plot. }
	\label{fig:dataPrivate_distancetooptimality}
\end{figure}

We start with discussing the data-hiding algorithm and take a look at its performance over iterations. Recall that most of the data exchange requirements are eliminated in this scheme. Parties need to share only their $\vs_k$ vectors at each iteration. In Figure \ref{fig:dataPrivate_distancetooptimality}, we report the percentage distances to the optimal objective function values of model \eqref{eqn:org_model_obj}-\eqref{eqn:org_model_nonneg} for the first 1,000 iterations. The solid line represents the mean whereas the shaded area shows the maximum and minimum values. Note that on the vertical axis of the figure, we present the percentage gaps with respect to the optimal objective function values.

As Figure \ref{fig:dataPrivate_distancetooptimality} shows, the objective function value obtained with the data-hiding algorithm can be less than $1\%$ away from the optimal objective function value. However, this percentage varies significantly as the problem instance changes. On average, when the number of iterations reaches $790$, the algorithm is already $25\%$ away from the optimal objective function value. This indicates that if parties use the algorithm step \eqref{eq:subup} only up to a certain number of iterations, they can still achieve an acceptable total utility by disclosing only their allotments from the shared capacities at each iteration. 

\begin{figure}[H]
	\centering
	\includegraphics[width=0.65\linewidth]{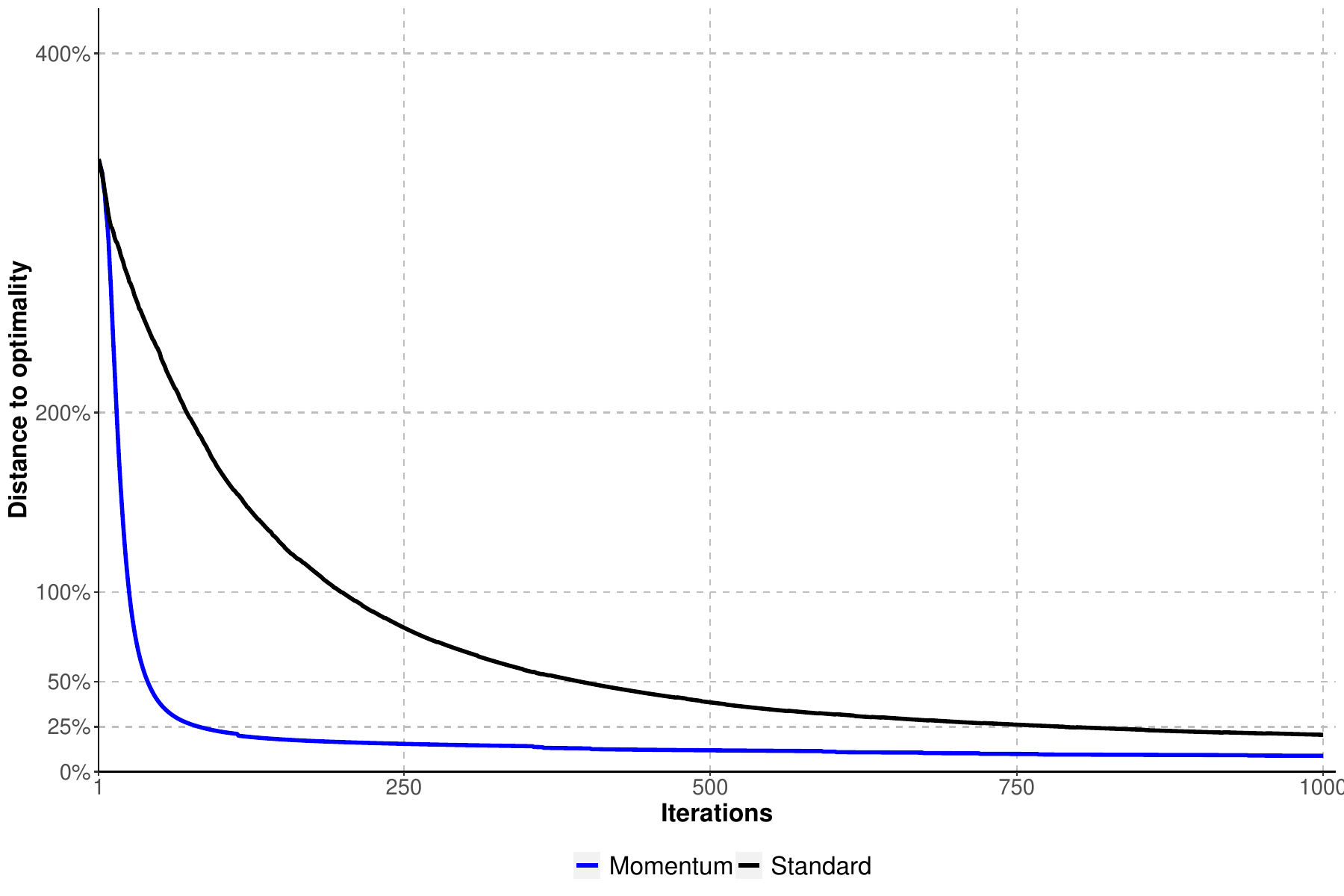}
	\caption{The comparison of data-hiding algorithms (with and without momentum updates) by showing their average percentage differences between the objective function values of the non-private model and the data-hiding models for the first 1,000 iterations.}
	\label{fig:comp_dataPrivate}
\end{figure}

We also include momentum updates on the data-hiding algorithm to investigate its impact on convergence. The results are given in Figure \ref{fig:comp_dataPrivate} for the first $1,000$ iterations. The figure depicts that momentum updates improve the convergence speed. On average, the algorithm with the momentum updates obtains objective function values that are $15\%$ away from the optimal one roughly in $284$ iterations. Without the momentum updates, the same level is not reached even after $1,000$ iterations. This shows the promise of using momentum updates with the data-hiding algorithm.

\begin{figure}[h]
	\centering
	\includegraphics[width=0.8\linewidth]{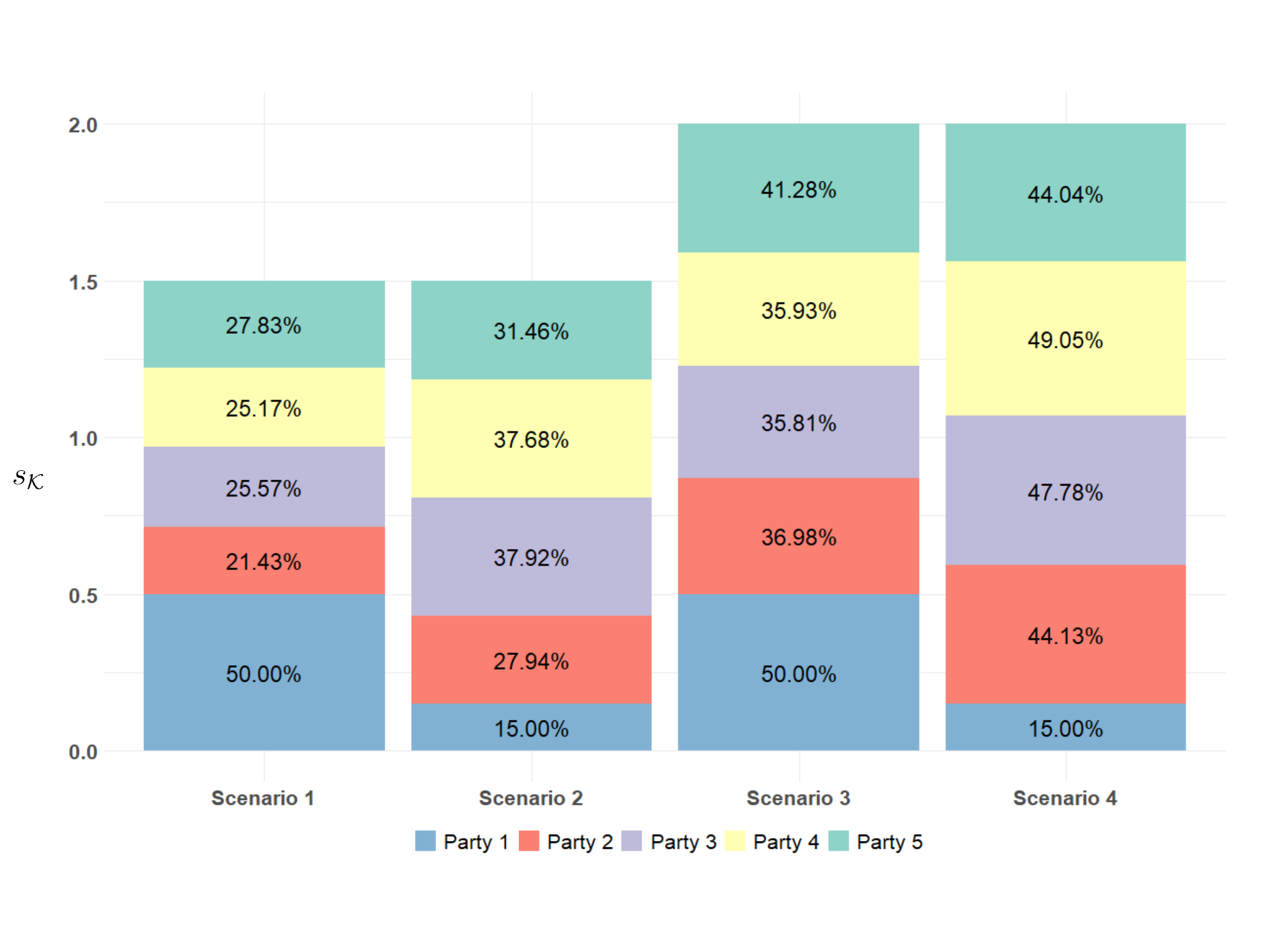}
	\caption{The shares of parties when $K = 5$.}
	\label{fig:party5_shares}
\end{figure}

\begin{figure}[t]
	\centering
	\includegraphics[width=0.8\linewidth]{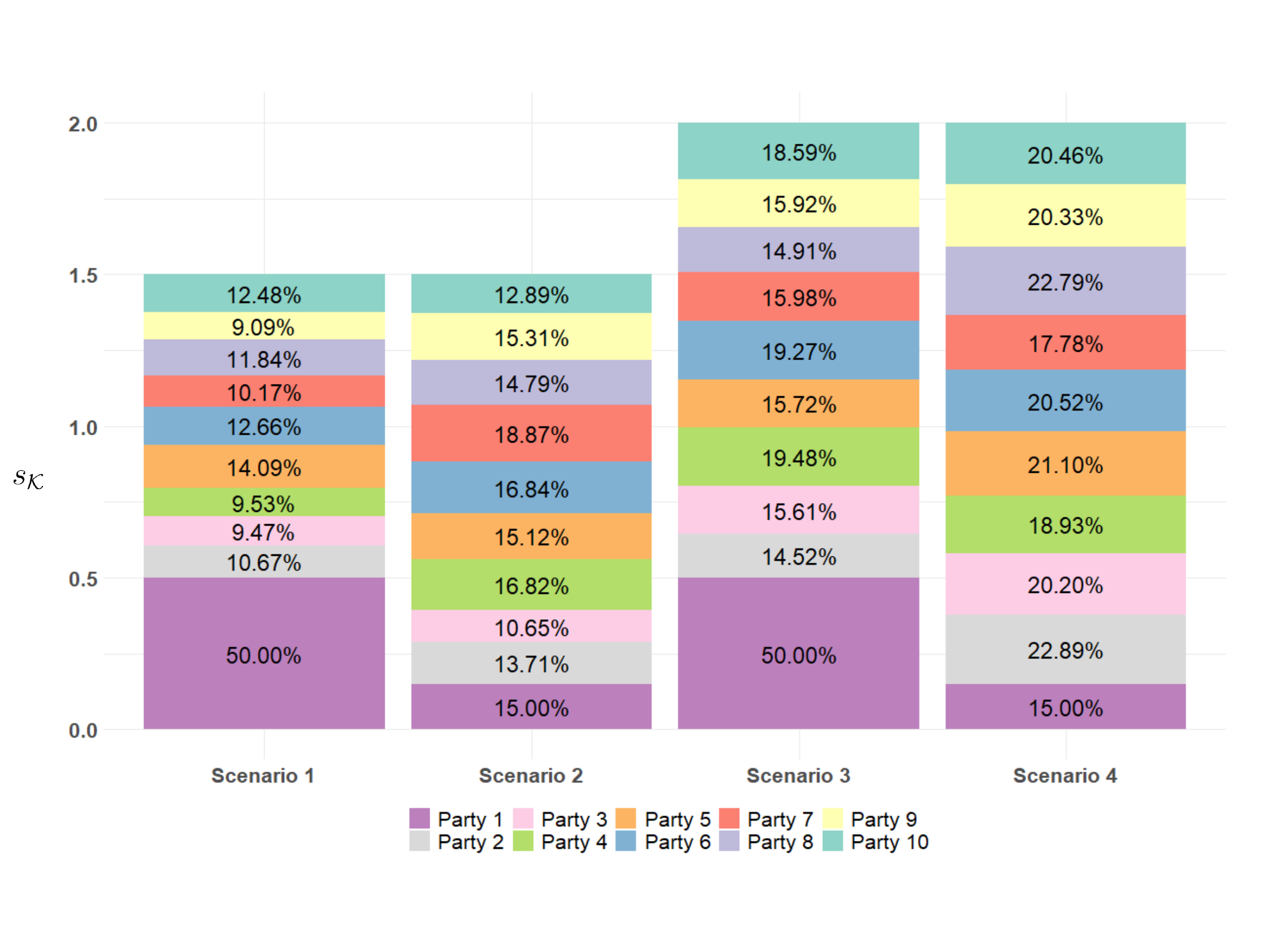}
	\caption{The shares of parties when $K = 10$.}
	\label{fig:party10_shares}
\end{figure}

Next, we conduct experiments with the differentially private algorithm and present our results. We analyze the effects of parameters ($\epsilon, K$) on the resulting objective function values and the optimality gap which is also known as the percentage regret in privacy literature \citep{Chen23}. To assess the impact of privacy on various market dynamics, we examine two scenarios. In the first scenario, we assume that a dominant participant, denoted as Party 1 $(k=1)$, is sufficiently influential to assert a claim of 50\% of the shared resources, represented as $\vs_{1} \leq 0.50 \vc$. The second scenario, where $\vs_{1} \leq 0.15 \vc$, simulates cases of moderate or low claims depending on the number of participating parties. The upper bounds for the market shares of the remaining parties 
are randomly generated. We also explore two scenarios to examine the impact of buffer resources utilized for compensating potential capacity overflow. For those scenarios, we consider an upper bound $\vs_{\mathcal{K}}$ for the total allocations for all parties, such that $\sum_{k = 1}^{K} \vs_{k} \leq \vs_{\mathcal{K}}$. In the first scenario, we analyze a moderate market where the total allocations for all parties can be at most 50\% above the levels of the shared resources; that is $\vs_{\mathcal{K}} = 1.5\vc$. Then, we analyze an acquisitive scenario where $\vs_{\mathcal{K}} = 2.0 \vc$. By undertaking this analysis, we aim to observe the influence of diverse buffer resource policies on the optimality gap. Figure \ref{fig:party5_shares} and Figure \ref{fig:party10_shares} illustrate these market scenarios across different problem sizes. 

\begin{figure}[H]
	\centering
	\includegraphics[width=0.95\linewidth]{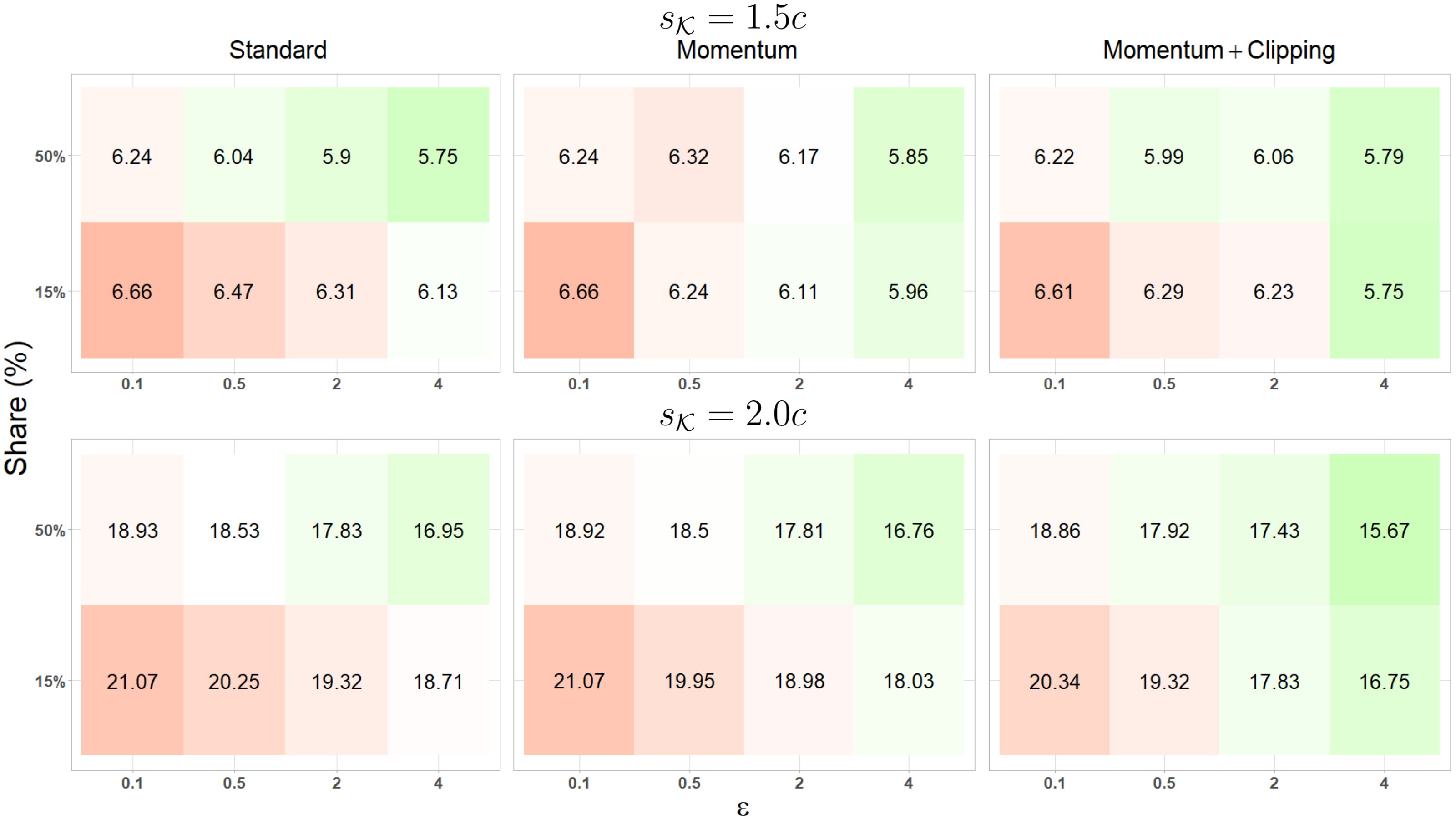}
	\caption{The average percentage differences between the objective function values of the non-private and the differentially private models for $K = 5$ and varying $\epsilon$ values.}
	\label{fig:party5}
\end{figure}

\begin{figure}[H]
	\centering
	\includegraphics[width=0.95\linewidth]{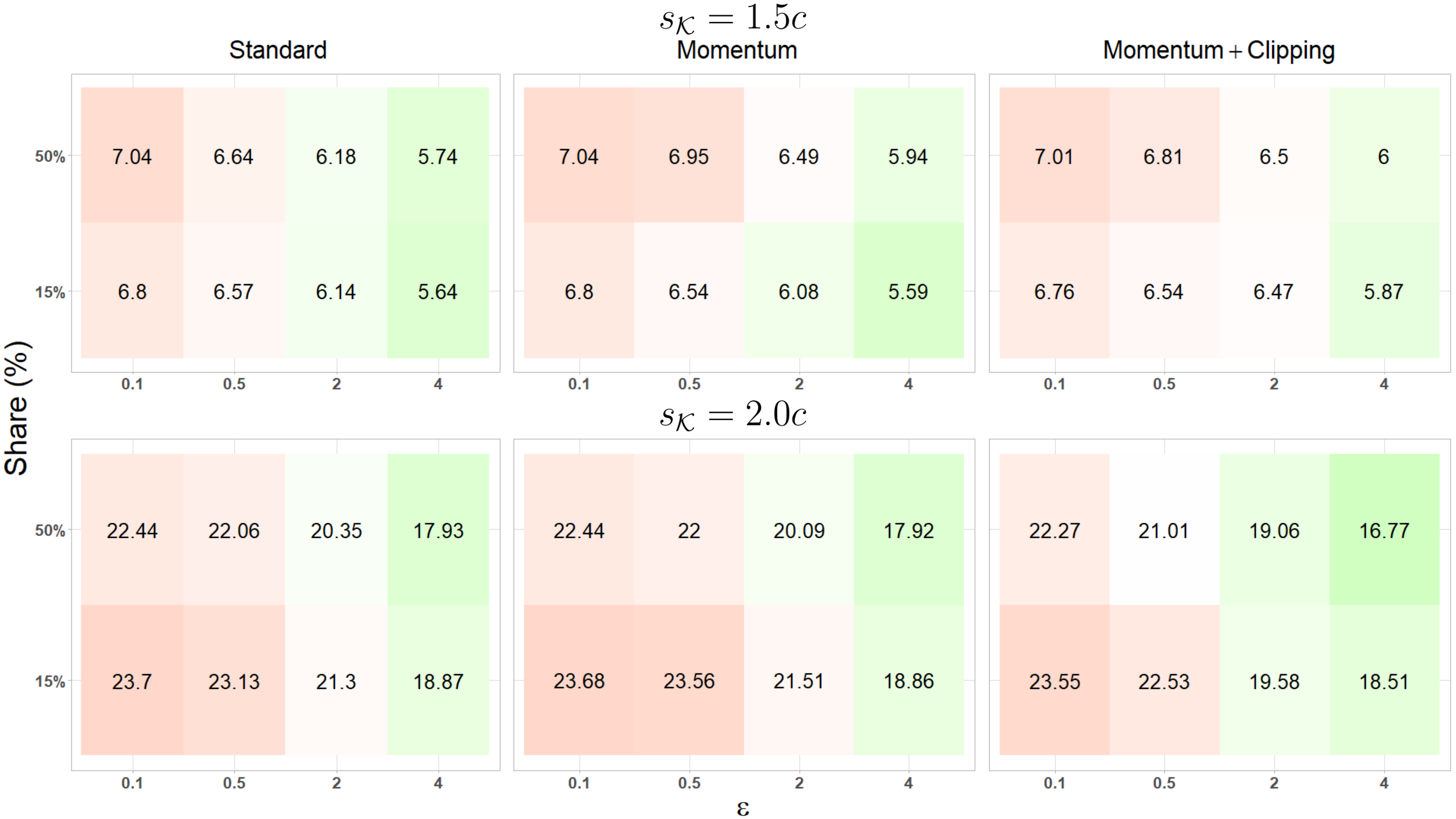}
	\caption{The average percentage differences between the objective function values of the non-private and the differentially private models for $K = 10$ and varying $\epsilon$ values.}
	\label{fig:party10}
\end{figure}

We set the maximum number of iterations to $T=150$ and the momentum parameter to $\gamma = 0.1$ for all scenarios and adjust the other parameters accordingly. Recall that momentum updates use information from past steps of the algorithm. As we add noise at each iteration of the algorithm, this noise also accumulates especially during the early iterations with the momentum updates. Moreover, the algorithm is run for only a limited number of iterations and the initial oscillations of the momentum updates cannot be dampened especially when the variance of the noise is high. Therefore, in these experiments, we set the momentum parameter to a low value. Figures \ref{fig:party5} and \ref{fig:party10} plot the average percentage gaps between the objective function values of the non-private and the differentially private models. 

In these figures, we investigate the performances of the differentially private models with standard and momentum updates by comparing the respective optimality gaps with the non-private model. We also present our results for the differentially private algorithm with momentum updates and clipping. The first observation we have is that the average percentage gaps between the non-private and the differentially private models decrease as the algorithm becomes less private. This phenomenon is attributed to the smaller magnitude of variances added on top of $\vs_k$ values. Parameter $\epsilon$ is known as the privacy budget and it is an important parameter for differential privacy protection. Altering $\epsilon$ allows us to explore the impact of privacy on the total utility. Smaller values of $\epsilon$ imply stronger privacy and a greater decrease in model utility, leading to a higher optimality gap or regret. Several studies investigate the effect of $\epsilon$ and point out that while a higher $\epsilon$ tends to yield improved model utility, it also corresponds to a reduced level of privacy protection \citep{Ryu22}. This trade-off necessitates careful consideration when designing differentially private multi-party resource-sharing mechanisms, as practitioners must strike a balance that aligns with their specific privacy and utility requirements. When we analyze the impact of buffer resources and the problem size, we obtain interesting findings. In the context of five parties, the market share scenarios $\vs_{1} \leq 0.50\vc$ and $\vs_{1} \leq 0.15\vc$ correspond to instances where Party 1 either dominates the market or is dominated by other parties, respectively. When a single party dominates the market, the average optimality gap or regret tends to be smaller. In a ten-party setup, the scenario $\vs_{1} \leq 0.15\vc$ corresponds to a case where none of the parties are dominant (see Figures \ref{fig:party5_shares} and \ref{fig:party10_shares}). Similarly, we observe a similar trend where the gap is generally smaller when one of the parties has a higher market share.

Next, we compare the results of standard differential privacy, momentum update, and momentum update with adaptive clipping. Figures \ref{fig:party5} and \ref{fig:party10} illustrate that their performances are generally similar. The impact of the momentum update becomes more pronounced with an increase in the number of parties and the utilization of buffer resources. Recall that different from the literature, our problem setting does not involve a trusted aggregator. Implementing differential privacy without a trusted aggregator leads to more substantial noise in the shared resources, consequently amplifying the optimality gap. This is more striking when the number of parties is high and the market exhibits an acquisitive nature ($\vs_{\mathcal{K}} = 2.0\vc$). This aligns with our theoretical results, as a larger dataset inherently leads to greater variance. As illustrated by Figures \ref{fig:party5} and \ref{fig:party10} clipping and adaptive sensitivity enhance the efficiency of the private decomposition algorithm by reducing the variance of the noise, and consequently, improve the optimality gaps, particularly when the privacy budget is low ($\epsilon \in \{ 0.1, 0.5\}$). The main advantage of this modified algorithm is the lower magnitude of variance added to $\vs_k$ values. In addition to the altered update steps, we truncate the $\tilde{\vs}_k$ values. This benefits from the post-processing property of a differentially private mechanism. Before sharing the noisy vector $\tilde{\vs}_k$, the party $k$ can truncate the noisy vector so that $\tilde{\vs}_k \in [\epsilon_{\vc}, \vc]$. This would not cause any privacy problems since differential privacy is immune to post-processing \citep{dwork2014algorithmic}.

\section{Conclusion} \label{sec:conclusion}

We have formulated a multi-party resource-sharing model, where the primary question is how to ensure the data privacy of the involved parties. We have started with a relaxation of the original model to obtain a decomposition that minimizes the necessity of sharing individual input data. This relaxation approach has led to a data-hiding scheme. Then, we have modified this scheme to benefit from the mathematically well-defined framework of differential privacy. This modification has allowed us to give a privacy guarantee and obtain our local differentially private resource-sharing algorithm. In local differential privacy, parties retain direct control over data privatization, removing the need for a trusted aggregator. This approach diverges from the standard differential privacy, where parties entrust their data to a central entity. Our work is positioned as a first step towards formally addressing data privacy without reliance on any trusted party in the context of multi-party resource sharing through optimization.

As differential privacy guarantee comes inevitably at the loss of optimality, we have also given bounds on this loss for different privacy parameters of the proposed algorithm. We have additionally proposed a modification to the original algorithm to increase the efficiency when the number of parties increases. Through a computational study on a planning problem with synthetic data, we have validated our analysis, presented results for various privacy parameters, and given a discussion on the impact of privacy levels on optimality. Our findings for the original algorithm have revealed that, as the number of parties decreases, the optimality gaps also diminish. Importantly, our proposed modification ensures that the optimality gap does not escalate significantly with an increasing number of parties.

In our view, the optimization perspective on differential privacy put forward in this paper opens up several opportunities for future research. We have reserved this work for a linear programming formulation because linear programs are used frequently in various applications. The methodology presented in this paper could be extended to a more general setting, in which the objective function and the constraint functions are possibly nonlinear but define a convex programming model. In that case, one needs to give a clear discussion about what constitutes private and non-private datasets. Combined with an application, this could be an interesting research topic. Any iterative differentially private algorithm loses privacy as the number of iterations increases. Therefore, our differentially private algorithm also stops after a fixed number of iterations. This early termination affects not only the optimality but also the feasibility of the algorithm. In other words, the resulting set of allocations from our algorithm may exceed some of the shared capacities. We are not aware of another approach that would lead to a differentially private algorithm that also guarantees feasibility. Investigating this point further is certainly on our agenda for future research.


\clearpage

\appendix

\section{Omitted Proofs} \label{sec:appProof}
We reserve this section for the omitted proofs in the main text.
\subsection{LDP Proof of Algorithm \ref{alg:dpra}} \label{sec: Proofs of the results related to LDP}

\propnum{\ref{prop:diffPriv}}{
The multi-party resource sharing algorithm using updates \eqref{eq:subupdp} for $T$ iterations provides $(\epsilon, \delta)$-LDP if the random vectors are drawn as 
\begin{equation*} 
\bm{\omega}_{k, j}^{(t)} \sim \mathcal{N}\left(0, \frac{T m  c_{i}^{2}}{2 \rho(\epsilon, \delta)} \right), \quad k = 1, \ldots, K; ~ j = 1, \ldots, m; ~t=1, \dots, T,
\end{equation*}
where $\rho(\epsilon, \delta) = (\sqrt{\log(1/\delta) + \epsilon} - \sqrt{\log(1/\delta)})^{2}$.}

\begin{proof}
To establish the $(\epsilon, \delta)$-LDP of the overall algorithm, we will show that each party $k$ shares its output in $(\epsilon, \delta)$-LDP during the whole course of the algorithm with $T$ iterations. 
Firstly, recall that the $L_{2}$-sensitivity of calculating $s_{k, j}^{(t)}$ is $c_{j}$ as indicated in \eqref{eq: bound on s}. Hence, by Theorem \ref{thm:gausszCDP}, for each component $j$, outputting $\hat{s}_{k, j}^{(t)} = s_{k, j}^{(t)} + w_{k, j}^{(t)}$ in \eqref{eq:subupdp} is $\rho_{0}$-zCDP where $\rho_{0} = \rho(\epsilon, \delta)/(T m)$. Since there are $m$ components in $s_{k}$, \eqref{eq:subupdp} becomes $m \rho_{0} =\rho(\epsilon, \delta)/T$-zCDP by Theorem  \ref{thm:zCDPComp}. Moreover, the update in \eqref{eq:subupdp} is performed for $T$ iterations, so the overall algorithm is  $T m \rho_{0} = \rho(\epsilon, \delta)$-zCMP. Finally, using the conversion in Theorem \ref{thm:zCDPtoLDP} and checking that $\rho(\epsilon, \delta) + \sqrt{4 \rho(\epsilon, \delta) \log(1/\delta)} = \epsilon$, we end up with $(\epsilon, \delta)$-LDP.
\end{proof}

\cornum{\ref{cor:difPrivClip}}{
The multi-party resource sharing algorithm using updates \eqref{eq:subupdpClip-1}-\eqref{eq:subupdpClip-4} for $T$ iterations is $(\epsilon, \delta)$-LDP.
}
\begin{proof}
Note that the sensitivity of $\min\{\bar{s}_{k, j}, \bar{s}_{k, j}\}$ in \eqref{eq:subupdpClip-2} is $\bar{s}_{k, j}$. Hence, by Theorem \ref{thm:gausszCDP}, for each component $j$, outputting $\hat{s}_{k, j}^{(t)} = s_{k, j}^{(t)} + w_{k, j}^{(t)}$ in \eqref{eq:subupdpClip-2} is $\rho_{0}$-zCDP where $\rho_{0} = \rho(\epsilon, \delta)/(T m)$. The rest of the proof follows the same lines as the proof of Proposition \ref{prop:diffPriv}.
\end{proof}

\subsection{Optimality Gap Proofs} \label{sec: Proofs of the theorems on the  optimality gap}
Next, we prove the theorems on the upper bounds of the optimality gap for the proposed methods. Before proceeding, we let $G^{(t)} =  \CL(\vx^{(t)}, \vs^{(t)}, \vlambda^{(t)}) - \CL(\vx^*, \vs^*, \vlambda^*)$ for $t \geq 0$ and $G^{(t_{1}, t_{2})} =  G^{(t_{1})} - G^{(t_{2})} = \CL(\vx^{(t_{1})}, \vs^{(t_{1})}, \vlambda^{(t_{1})}) - \CL(\vx^{(t_{2})}, \vs^{(t_{2})}, \vlambda^{(t_{2})})$ for $t_{2} \geq t_{1} \geq 0$.


\thmnum{\ref{thm:convergence}}{
Suppose Assumption \ref{asmp: bound on lambda} holds and Algorithm \ref{alg:dpra} is run for $T$ iterations with a fixed step-size $\nu^{(t)} = \nu$ and with noise variances given in \eqref{eq: capacity-specific DP variances-1} to provide $(\epsilon, \delta)$-LDP. Then, there exists a $\nu$, for which the minimum expected distance to the optimal objective function value is bounded as
\[
\min_{t=0, 1, \ldots, T}~ \expec[G^{(t)}] \leq \frac{M B}{\sqrt{T}},
\]
where the constant $B$ is defined as
\begin{equation*} 
B = \|\vc\| \left(\frac{T m K}{2\rho(\epsilon, \delta)}+(K-1)^2 \right)^{1/2}.
\end{equation*}}

\begin{proof}
Let $\rho = \rho(\epsilon, \delta)$. Recall the differentially private update steps \eqref{eq:subupdp} and $\tilde{\vs}_k^{(t)} = \vs_k^{(t)} + \bm{\omega}_k^{(t)}$, where $\omega_{k, j}^{(t)} \sim \mathcal{N}(0, Tm c_{j}^{2}/(2\rho))$ independently for $j = 1,\ldots, m$. Let 
\begin{equation} \label{eq: approximate subgradient}
\tilde{\vg}^{(t)} = \vc - \sum_{k = 1}^{K} \tilde{\vs}_k^{(t)}
\end{equation}
denote the approximate subgradient. Then, we have 
\begin{equation} \label{eq:expecNoisy}
\begin{aligned}
\expec[\|\tilde{\vg}^{(t)}\|^2] &\leq\frac{Tm K \|\vc\|^2}{2\rho} + \big\|(K-1)\vc \big\|^2 \\
&= \left(\frac{Tm K}{2\rho}+(K-1)^2\right) \|\vc\|^2\\
&= B^{2},
\end{aligned}
\end{equation}
where $B^{2}$ was defined in \eqref{eq: B bound}. We also have
\begin{align*}
\expec[\|\vlambda^{(t+1)} - \vlambda^*\|^2 | \vlambda^{(t)}] &= \expec[\|\vlambda^{(t)} - \vlambda^* - \nu^{(t)} \tilde{\vg}^{(t)}\|^2 | \vlambda^{(t)}] \\[1em] & =\|\vlambda^{(t)} - \vlambda^*\|^2 - 2\nu^{(t)}\expec[\tilde{\vg}^{(t)}| \vlambda^{(t)}]\tr(\vlambda^{(t)} - \vlambda^*) + (\nu^{(t)})^2\expec[\|\tilde{\vg}^{(t)}\|^2 | \vlambda^{(t)}].
	\end{align*}
	Using the definition of subgradient, we have 
	\[\CL(\vx^{*}, \vs^{*}, \vlambda^{*}) \geq \expec[\tilde{\vg}^{(t)}| \vlambda^{(t)}]\tr(\vlambda^{*} - \vlambda^{(t)}) + \CL(\vx^{(t)}, \vs^{(t)}, \vlambda^{(t)}),
	\]
	which yields
	\begin{align}
		\expec[\|\vlambda^{(t+1)} - \vlambda^*\|^2 | \vlambda^{(t)}] \leq\|\vlambda^{(t)} - \vlambda^*\|^2 - 2\nu^{(t)} G^{(t)} + (\nu^{(t)})^2\expec[\|\tilde{\vg}^{(t)}\|^2| \vlambda^{(t)}]. \label{eqn:lastLine}
	\end{align}
Then, taking the expectation of both sides of \eqref{eqn:lastLine} leads to
\begin{equation}
\label{eqn:beforeB}	
2\nu^{(t)}\expec[G^{(t)}] \leq \expec[\|\vlambda^{(t)} - \vlambda^*\|^2] - \expec[\|\vlambda^{(t+1)} - \vlambda^*\|^2] + (\nu^{(t)})^2\expec[\|\tilde{\vg}^{(t)}\|^2].
\end{equation}
By \eqref{eq:expecNoisy} we use $B^2$ for an upper bound for the term $\expec[\|\tilde{\vg}^{(t)}\|^2]$ in \eqref{eqn:beforeB}, and take the summation of both sides to obtain
	\[
	2\sum_{t=0}^{T-1}\nu^{(t)}\expec[G^{(t)}] \leq \expec[\|\vlambda^{(0)} - \vlambda^*\|^2] - \expec[\|\vlambda^{(T)} - \vlambda^*\|^2] +B^2 \sum_{t=0}^{T-1}(\nu^{(t)})^2 \leq \expec[\|\vlambda^{(0)} - \vlambda^*\|^2] + B^2 \sum_{t=0}^{T-1}(\nu^{(t)})^2.
	\]
	This inequality also implies
	\begin{equation}
	\label{eqn:finalBound}
	\min_{t=0, \dots, T-1} \expec[G^{(t)}] \leq \frac{\expec[\|\vlambda^{(0)} - \vlambda^*\|^2] +B^2 \sum_{t=0}^{T-1}(\nu^{(t)})^2}{2 \sum_{t=0}^{T-1}\nu^{(t)}}.
	\end{equation}
	Recall that our assumption satisfies $\|\vlambda^{(0)} - \vlambda^*\|^2\leq M^2$. With the constant step-size $\nu^{(t)} = \frac{M}{B\sqrt{T}}$, inequality \eqref{eqn:finalBound} becomes
	\[
	\min_{t=0, \dots, T-1} \expec[G^{(t)}]  \leq \frac{M^2 +B^2 \sum_{t=0}^{T-1}(\nu^{(t)})^2}{2\sum_{t=0}^{T-1}\nu^{(t)}} = \frac{M B}{\sqrt{T}} = M\|\vc\| \sqrt{\frac{mK}{2\rho} + \frac{(K-1)^2}{T}}.
	\]
This shows the desired result.
\end{proof}

\thmnum{\ref{thm:convergence_tight}}{
Suppose Assumption \ref{asmp: bound on lambda} holds. If the multi-party resource-sharing algorithm with momentum updates is run with a fixed step-size for $T$ iterations with $(\epsilon, \delta)$-LDP, the minimum expected distance to the optimal objective function value is bounded as
\[
\min_{t=0, 1, \ldots, T}~ \expec [G^{(t)}] \leq \frac{M^2}{2 T \nu} + \frac{\nu}{2} B^{2},
\]
provided that the step-size satisfies $0 \leq \nu \leq \frac{\sqrt{ G^{(0)2} + (1 - \gamma) B^{2}T M^{2}} - G^{(0)}}{B^{2} T}$ with $B$ as in \eqref{eq: B bound}.
}


\begin{proof}
The update step presented in \eqref{eqn:dp_momentum} can be rewritten as 
\[
\vlambda^{(t+1)} = \vlambda^{(t)} - \nu^{(t)} \tilde{\vg}^{(t)} + \momentumParam \left(\vlambda^{(t)} - \vlambda^{(t-1)}\right),
\]
where the approximate subgradient $\tilde{\vg}^{(t)}$ was defined in \eqref{eq: approximate subgradient}.
We use an equivalent representation of the update step in \citep{ghadimi2015global} and obtain
\[
\vlambda^{(t+1)} + \vp^{(t+1)} = \vlambda^{(t)} + \vp^{(t)} - \frac{\nu^{(t)}}{1 - \momentumParam} \tilde{\vg}^{(t)},
\]
where
\[
\vp^{(t)}=\begin{cases} 
\frac{\momentumParam}{1-\momentumParam}\left(\vlambda^{(t)}-\vlambda^{(t-1)}\right), & t \geq 1; \\
0, & t = 0.
\end{cases}
\]
Then, we proceed as
\begin{align*}
	& \expec\left[\|\vlambda^{(t+1)} + \vp^{(t+1)} -\vlambda^*\|^2 | \vlambda^{(t)}\right] = \expec\left[\|\vlambda^{(t)} + \vp^{(t)} -\vlambda^* - \frac{\nu^{(t)}}{1-\momentumParam}\tilde{\vg}^{(t)}\|^2| \vlambda^{(t)}\right] \\
	& = \|\vlambda^{(t)} + \vp^{(t)} -\vlambda^*\|^2 - \expec\left[\frac{2\nu^{(t)}}{1-\momentumParam}(\vlambda^{(t)} + \vp^{(t)} -\vlambda^*)\tr\tilde{\vg}^{(t)} - \left(\frac{\nu^{(t)}}{1-\momentumParam}\right)^2 \|\tilde{\vg}^{(t)}\|^2 | \vlambda^{(t)}\right]\\
	& = \|\vlambda^{(t)} + \vp^{(t)} -\vlambda^*\|^2 -  \expec\left[\frac{2\nu^{(t)}}{1-\momentumParam}(\vlambda^{(t)} - \vlambda^*)\tr\tilde{\vg}^{(t)} - \frac{2\nu^{(t)}}{1-\momentumParam}\left(\frac{\momentumParam}{1-\momentumParam}(\vlambda^{(t)}-\vlambda^{(t-1)})\right)\tr \tilde{\vg}^{(t)} | \vlambda^{(t)}\right]  \nonumber\\ & \qquad\qquad + \expec\left[\left(\frac{\nu^{(t)}}{1-\momentumParam}\right)^2 \|\tilde{\vg}^{(t)}\|^2 | \vlambda^{(t)} \right] \\
	& = \|\vlambda^{(t)} + \vp^{(t)} -\vlambda^*\|^2 - \frac{2\nu^{(t)}}{1-\momentumParam} \expec\biggl[ (\vlambda^{(t)} - \vlambda^*)\tr\tilde{\vg}^{(t)}  | \vlambda^{(t)} \biggr] - \frac{2\nu^{(t)}\momentumParam}{(1-\momentumParam)^2} \mathbb{E} \left[ \left(\vlambda^{(t)}-\vlambda^{(t-1)}\right)\tr \tilde{\vg}^{(t)} | \vlambda^{(t)}\right] \\
& \qquad \qquad + \left(\frac{\nu^{(t)}}{1-\momentumParam}\right)^2 \mathbb{E} \left[\|\tilde{\vg}^{(t)}\|^2 | \vlambda^{(t)}\right].
\end{align*}
Using the definition of subgradient, we obtain 
\[\CL(\vx^{*}, \vs^{*}, \vlambda^{*}) \geq \expec[\tilde{\vg}^{(t)}| \vlambda^{(t)}]\tr(\vlambda^{*} - \vlambda^{(t)}) + \CL(\vx^{(t)}, \vs^{(t)}, \vlambda^{(t)}),
\]
which yields
\begin{align*}
 \expec\left[\|\vlambda^{(t+1)} + \vp^{(t+1)} -\vlambda^*\|^2 | \vlambda^{(t)}\right] & \leq \|\vlambda^{(t)} + \vp^{(t)} -\vlambda^*\|^2 -  \frac{2\nu^{(t)}}{1-\momentumParam} G^{(t)}  \nonumber\\
 & \qquad + \expec\biggl[-   \frac{2\nu^{(t)}\momentumParam}{(1-\momentumParam)^2}\left(\vlambda^{(t)}-\vlambda^{(t-1)}\right)\tr \tilde{\vg}^{(t)} + \left(\frac{\nu^{(t)}}{1-\momentumParam}\right)^2 \|\tilde{\vg}^{(t)}\|^2 | \vlambda^{(t)}\biggr].
\end{align*}
Again from the definition of subgradient, we have
\[\CL(\vx^{(t-1)}, \vs^{(t-1)}, \vlambda^{(t-1)}) \geq \expec[\tilde{\vg}^{(t)}| \vlambda^{(t)}]\tr(\vlambda^{(t-1)} - \vlambda^{(t)}) + \CL(\vx^{(t)}, \vs^{(t)}, \vlambda^{(t)}).
\]
We end up with
\begin{align*}
	\expec\left[\|\vlambda^{(t+1)} + \vp^{(t+1)} -\vlambda^*\|^2 | \vlambda^{(t)}\right] & \leq \|\vlambda^{(t)} + \vp^{(t)} -\vlambda^*\|^2  -  \frac{2\nu^{(t)}}{1-\momentumParam} G^{(t)} + \frac{2\nu^{(t)}\momentumParam}{(1-\momentumParam)^2} G^{(t-1, t)}  \nonumber\\ 
& \qquad \qquad + \left(\frac{\nu^{(t)}}{1-\momentumParam}\right)^2  \expec\left[ \|\tilde{\vg}^{(t)}\|^2 | \vlambda^{(t)}\right].
\end{align*}
We then take the expectation of both sides
\begin{align*}
	\frac{2\nu^{(t)}}{1-\momentumParam}\expec\left[G^{(t)} \right] - \frac{2\nu^{(t)}\momentumParam}{(1-\momentumParam)^2}\expec\left[ G^{(t-1, t)} \right] & \leq \expec\left[\|\vlambda^{(t)} + \vp^{(t)} -\vlambda^*\|^2\right] - \expec\left[\|\vlambda^{(t+1)} + \vp^{(t+1)} -\vlambda^*\|^2 \right] \nonumber\\ 
	& + \left(\frac{\nu^{(t)}}{1-\momentumParam}\right)^2 \expec\left[\|\tilde{\vg}^{(t)}\|^2\right].
\end{align*}
If we take the step-size parameter fixed, $\nu^{(t)} = \nu$, $G^{(t)} = G^{(0)}$ for $t \leq 0$, and sum over iterations $t = 0, \ldots, T-1$, we obtain
\begin{align*}
	\frac{2 \nu}{1-\momentumParam} \sum_{t=0}^{T-1} \expec\left[ G^{(t)} \right] - \frac{2\nu\momentumParam}{(1-\momentumParam)^2}\expec\left[ G^{(0, T-1)} \right] & \leq \expec\left[\|\vlambda^{(0)} -\vlambda^*\|^2\right] - \expec\left[\|\vlambda^{(T)} + \vp^{(T)} -\vlambda^*\|^2 \right] +  T \left(\frac{\nu}{1-\momentumParam}\right)^2 B^{2} \nonumber \\
	& \leq \expec\left[\|\vlambda^{(0)} -\vlambda^*\|^2\right]+ T\left(\frac{B\nu}{1-\momentumParam}\right)^2,
\end{align*}
since $\expec[\|\tilde{\vg}^{(t)}\|^2] \leq B^2$ by \eqref{eq:expecNoisy}. This is equivalent to the following
\begin{align*}
	2 \nu\sum_{t=0}^{T-1}\expec\left[G^{(t)} \right] \leq (1-\momentumParam)\expec\left[\|\vlambda^{(0)} -\vlambda^*\|^2\right] + \frac{2\nu\momentumParam}{(1-\momentumParam)}\expec\left[ G^{(0, T-1)}\right] + \frac{T(\nu B)^2}{1-\momentumParam}.
\end{align*}
This inequality also implies
\begin{align*}
\min_{t=0, \dots, T-1} \expec[ G^{(t)}] \leq \frac{(1-\momentumParam)\expec\|\vlambda^{(0)} - \vlambda^*\|^2 +\frac{2\nu\momentumParam}{(1-\momentumParam)}\expec\left[ G^{(0, T-1)} \right] + \frac{T(\nu B)^2}{1-\momentumParam}}{2T\nu}.
\end{align*}
Recall the assumption $\|\vlambda^{(0)} - \vlambda^*\|^2\leq M^2$, which implies that $G^{(0)} < \infty$. This further implies that $\mathbb{E}[G^{(0, t)}] \leq G^{(0)}$. Therefore, 
\begin{align*}
	\min_{t=0, \dots, T-1}  \mathbb{E}[G^{(t)}] \leq \frac{(1-\momentumParam)M^2 +\frac{2\nu\momentumParam}{1-\momentumParam}G^{(0)} + \frac{T(\nu B)^2}{1-\momentumParam}}{2T\nu}.
\end{align*}
To obtain a tighter bound, we need to satisfy the following inequality
 \begin{align*}
	\frac{(1-\momentumParam)M^2 +\frac{2\nu\momentumParam}{1-\momentumParam} G^{(0)} + \frac{T(\nu B)^2}{1-\momentumParam}}{2T\nu} \leq \frac{M^2 +B^2 T\nu^2}{2T\nu}. 
\end{align*}
Thus, if we select $\nu$ as
\[
0 \leq \nu \leq \frac{\sqrt{ G^{(0)2} + (1 - \gamma) B^{2}T M^{2}} - G^{(0)}}{B^{2} T}, 
\]
then we guarantee to have a tighter bound.
\end{proof}

\section{Dual Formulation}
\label{sec:appDual}
Let $\valpha_1, \ldots, \valpha_{K}$ and $\vbeta_{1}, \ldots, \vbeta_{K}$ be the dual vectors corresponding to the constraints \eqref{eqn:sec_model_Ak} and \eqref{eqn:sec_model_nonnegprev}, respectively. Further, we define $\vlambda$ as the dual vector for Constraint \eqref{eqn:sec_model_sk_c}. Then, the dual formulation of \eqref{eqn:sec_model_obj}-\eqref{eqn:sec_model_nonneg} becomes 
\begin{align*}
\minimize & \vc\tr\vlambda + \sum_{k = 1}^{K} \vb_k\tr\vbeta_k \\
\subto &  \mA_k\tr \valpha_k + \mB_k\tr\vbeta_k \geq \vu_k, & k = 1, \ldots, K, \\
& \vlambda - \valpha_k \geq \zv, & k = 1, \ldots, K, \\
& \valpha_k, \vbeta_k \geq \zv, & k = 1, \ldots, K, \\
& \vlambda \text{ unrestricted}.
\end{align*}

\clearpage

{\footnotesize
\singlespace
\bibliographystyle{apalike}
\bibliography{./bib/privacyRM}
}

\end{document}